\documentclass[11pt,openany]{article} 
\usepackage{amssymb} \usepackage{amsmath,amscd,mathrsfs}

\usepackage{graphpap}
\usepackage{pstricks}
\usepackage{pst-all}
\usepackage{pstcol}
\usepackage{color}
\usepackage{enumitem}
\usepackage{hyperref}
\usepackage[active]{srcltx}

\newtheorem{thm}{Theorem}
\newtheorem{defin}{Definition}
 \newtheorem{cor}{Corollary}
\newtheorem{lem}{Lemma}
\newtheorem{rem}{Remark}
\newtheorem{prop}{Proposition}
\newtheorem{exam}{Example}

\def\N{{ \! \rm \ I\!N}}

\def\R{{ \! \rm \ I\!R}}

   \def \square{\hbox {$\sqcup
    $\llap {$\sqcap $}}} 
  
\newcommand{\parcial}[2]{\frac{\partial#1}{\partial#2}}
\newcommand{\norm}[1]{\vert \vert #1 \vert \vert}
\newcommand{\normv}[1]{ \mid #1 \mid }

\title{Time dependent quantum scattering theory on
 complete manifolds with a corner of codimension
2}
\author{Leonardo A. Cano Garc\'{i}a}
\begin{document}
\maketitle
\begin{abstract}
We show the existence and orthogonality of wave operators
naturally associated to a compatible Laplacian on a complete
manifold with a corner of codimension 2. In fact,  we prove
asymptotic completeness i.e. that the image of these wave
operators is equal to the space of absolutely continuous states of
the compatible Laplacian. We achieve this last result using time
dependent methods coming from many-body Schr\"odinger equations.
\end{abstract}
\section{Introduction}\label{section}
In this article we use analytic tools to tackle problems of
quantum scattering theory naturally associated to  geometric
Laplacians, at the same time this makes explicit the interactions
between the geometry of the manifold and the quantum dynamics of
the Laplacians.

Classical mechanics  tells us that the time-asymptotic behavior of
$n$-particles interacting with a pairwise potential of short range
can be described by clusters whose centers of mass   do not "feel"
each other. In the papers \cite{SIGALSOFFER1} and
\cite{SIGALSOFFER2} it was proved that a similar phenomenon occurs
in  quantum mechanics for many-particle Schr\"odinger operators
with short range potentials. These proofs were time-dependent and
geometric in nature, and they were initially developed in the
papers~\cite{GRAF} and~\cite{Y1}. In this article we prove
asymptotic completeness for compatible Laplacians on complete
manifolds with corners of codimension 2, which we abbreviate
c.m.w.c.2 through the text, by adapting the proof of~\cite{Y1} as
explained in~\cite{HS1}. Even though the ideas are adapted in a
quite direct way, we believe that this article provides a deeper
understanding of the spectral theory of compatible Laplacians on
c.m.w.c.2 and of the geometric insight behind the proof of the
results in~\cite{GRAF} and~\cite{Y1}, since the spectral analysis
of Schr\"odinger operators and geometric Laplacians are analogous
but not exactly the same.

The motivation to study these manifolds is the same as
in~\cite{CANOTHESIS} and~\cite{CanoMourre}: they work as  toy
models for understanding singularities as those that appear on
symmetric spaces of rank greater than $0$; they are  natural
examples of complete manifolds whose spectral theory is well
known, since they are a natural geometric generalization of the
Cartesian products of complete manifolds with cylindrical ends.
This last class of manifolds is very important in the study of the
index theorems of the seminal paper~\cite{APS} and we believe that
a deeper understanding of the spectral theory  of compatible
Laplacians on c.m.w.c.2 (see section~\ref{sec:compatible Lapla on
man corner 2})
 will shed light on the nature of the generalization of such theorems, specifically in order
to complete the method applied in~\cite{MuellerCorner}.
Generalizations of the index theorems of~\cite{APS} to c.m.w.c.2
were obtained in~\cite{L2signMazzeo} using surgery methods, we
believe that these formulas are related to  our scattering
operator (see~(\ref{eq:defin scattering})). Finally, our work
shows a clear analogy between many--particle Schr\"odinger
operators and the compatible Laplacians on c.m.w.c.2, this analogy
provides  a deeper understanding of the geometric nature of the
spectral theory of the former operators.
\subsection{Compatible Laplacians on complete manifolds with a corner
of codimension 2} \label{sec:compatible Lapla on man corner
2}Following~\cite{MuellerCorner}, we explain the notions of {\it
compact and complete manifolds with a corner of codimension 2} as
is done in~\cite{CANOTHESIS} and~\cite{CanoMourre}. Let $X_0$ be a
compact oriented Riemannian manifold with boundary $M$ and suppose
that there exists a hypersurface $Y$ of $M$ that divides $M$ in
two manifolds with boundary $M_1$ and $M_2$, i.e. $M=M_1 \cup M_2$
  and $Y=M_1 \cap M_2$. Assume also that  a neighborhood of $Y$ in $M$ is diffeomorphic to $Y \times (-\varepsilon,\varepsilon)$. We say that the manifold $X_0$ {\bf has a corner of
codimension~2} if $X_0$ is endowed with a Riemannian metric $g$
that is a
  product metric on small neighborhoods, $M_i \times (-\varepsilon,0]$ of the $M_i$'s and
  on a small neighborhood $Y \times (-\varepsilon,0]^2$ of  the corner $Y$. If $X_0$ has a corner of codimension 2, we say that  $X_0$ is a {\bf compact manifold with a corner of codimension 2}.
\bigskip
\begin{center}
\psset{unit=0.5cm}
    \begin{pspicture}(4,1)(12,7)
            \pscurve[](6.5,4)(7,3.7)(7.5,3.5)(8,4)
                \pscurve[](7,3.7)(7.5,3.9)(7.5,3.5)
            \pscurve[](10,6)(7,5.8)(5,5.7)(7,5.5)(9,5.4)
            \pscurve[](10,6)(9.5,3)(9,5.4)
            \pscurve[](5,5.7)(6,2.8)(8.5,3)(9.5,3)
            \rput(7.2,6.5){$M_1$}
            \rput(10.8,4.6){$M_2$}
            \pscircle*(10,6){0.2}
            \pscircle*(9,5.4){0.2}
            \psline(9,5.4)(11.5,6.1)
            \psline(10,6)(11.5,6.1)
            \rput(11.9,6.3){$Y$}
            \rput(8,1.8){Figure 1. Compact manifold with a corner of codimension 2.}
        \end{pspicture}
\end{center}
\begin{exam} For $i=1,2$, let $M_i$ be a compact
oriented Riemannian manifold with boundary $\partial M_i:=Y_i$.
Suppose that on a neighborhood $Y_i \times (-\varepsilon, 0]$ of
$Y_i$ the Riemannian metric $g_i$ of $M_i$ is a product metric
i.e.
 $g_i:=g_{Y_i}+du \otimes du$ where $u$ is the coordinate
 associated to the interval $(-\varepsilon, 0]$ in
 $Y_i \times (-\varepsilon, 0]$ and $g_{Y_i}$ is a Riemannian
 metric on $Y_i$ independent of $u$. Then the Cartesian product $M_1 \times M_2$ is a compact manifold with a corner of codimension 2.\end{exam}
Throughout this article we will denote $\R_+:=[0,\infty)$. From
the compact manifold with a corner $X_0$ we construct a complete
manifold $X$. Let $ Z_i:=M_i \cup_Y (\R_+ \times Y), \text{
i=1,2}$, where the bottom $\{0\} \times Y$ of the half-cylinder
$\R_+ \times Y$ is identified with $\partial M_i=Y$. Then $Z_i$ is
a complete
 manifold with cylindrical end.
Let us define the manifolds
\begin{equation*}
W_1:=X_0\cup_{M_2} (\R_+ \times M_2) \text{ and
}W_2:=X_0\cup_{M_1} (\R_+ \times M_1).
\end{equation*}
Observe that $W_i$ is an $n$-dimensional manifold with boundary
$Z_i$ that can be equipped with a Riemannian metric compatible
with the product Riemannian metric of $\R_+ \times M_2$ and the
Riemannian metric of $X_0$. Let:
\begin{equation*}
X:=W_1 \cup_{Z_1}(\R_+ \times Z_1)=W_2 \cup_{Z_2}(\R_+ \times
Z_2),
\end{equation*}
where we  identify $\{0\}\times Z_i$ with $Z_i$, the boundary of
$W_i$.
\bigskip
\begin{center}
      \psset{unit=0.5cm}
          \begin{pspicture}(0,-1.5)(10,11)
                  \psline(0,0)(0,10)
                   \qline(0,0)(10,0)
        \psline(2,0)(2,10)
        \psline(0,2)(10,2)
        \rput(1,1){$X_0$}
        \rput(9,9){$[0,\infty)^2 \times Y$}
        \rput(9,2.8){$Z_1$ }
        \rput(2.6,10.5){$Z_2$}
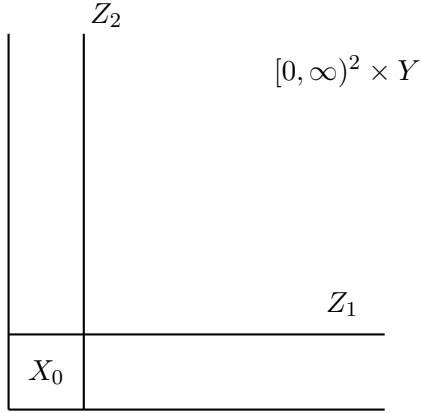
        \rput(6,-1){Figure 2. Sketch of a complete manifold with a corner of codimension 2.}
    \end{pspicture}
\end{center}
The  picture above is a sketch, in particular the lines that
enclose the picture should not be thought as boundaries.

Let $T \geq 0$ be given and set $Z_{i,T}:=M_i \cup_Y([0,T]\times
Y), \text{ for }i=1,2$, where $\{0\} \times Y$ is identified with
$Y$, the boundary of $M_i$. $Z_{i,T}$ is a family of manifolds
with boundary which exhausts $Z_i$. Next we attach to $X_0$ the
manifold $[0,T] \times M_1$ by identifying $\{0\} \times M_1$ with
$M_1$. The resulting manifold $W_{2,T}$ is a compact manifold with
a corner of codimension 2, whose boundary is the union of $M_1$
and $Z_{2,T}$. The manifold $X$ has associated a natural
exhaustion given by
\begin{equation}\label{eq: exhaustion X}
X_T:=W_{2,T} \cup_{Z_{2,T}}([0,T] \times Z_{2,T}), \text{ }T\geq
 0, \end{equation}
where we identify $Z_{2,T}$ with $\{0\} \times Z_{2,T}$.
\bigskip
\begin{center} \psset{unit=0.5cm}
    \begin{pspicture}(0,0)(24,10.5)
            \pscurve[](6.5,4)(7,3.7)(7.5,3.5)(8,4)
            \pscurve[](7,3.7)(7.5,3.9)(7.5,3.5)
        \pscurve[](5,5.7)(6,2.8)(8.5,3)(9.5,3)
        \pscurve[](10,11)(7,10.8)(5,10.7)(7,10.5)(9,10.4)
        \psline(9,10.4)(14,10.4)
        \psline(14,10.4)(14,5.3)
        \psline(9.5,3)(14.5,3)
        \psline(5,5.7)(5,10.7)
        \psline(10,11)(15,11)
        \psline(15,11)(15,6)
        \pscurve[](15,6)(14.5,3)(14,5.4)
        \rput(9,1){Figure 2. $X_T$, element of the exhaustion of $X$.}
        \psline(9,5.4)(9,10.4)
        \psline(5,5.7)(5,10.7)
        \psline[linestyle=dashed](10,6)(10,11)
        \rput(6.8,8.5){{\tiny $[0,T] \times M_1$}}
        \pscurve[linestyle=dashed](10,6)(8,5.9)(6,5.8)(5.5,5.9)(5,5.7)
\pscurve[](5,5.7)(6,5.4)(7,5.5)(9,5.4)
        \pscurve[linestyle=dashed](10,6)(9.8,3.4)(9.5,3)
\pscurve[](9.5,3)(9.2,3.2)(9,5.4)
        \pscurve[](5,5.7)(6,2.8)(8.5,3)(9.5,3)
        \pscurve[](15,6)(14.5,3)(14,5.4)
        \psline(9.5,3)(14.5,3)
        \psline[linestyle=dashed](15,6)(10,6)
        \psline(14,5.4)(9,5.4)
        \rput(11.5,4.7){{\tiny$ [0,T] \times M_2$}}
        \rput(12,8.5){{\tiny $ [0,T]^2\times Y$}}
        \rput(6,4.5){{\tiny $X_0$}}
\end{pspicture}
\end{center}
For each $T \in [0,\infty)$, $X$ has two submanifolds with
cylindrical ends, namely $(\{T\} \times M_i)   \cup (\{T\}\times
[0,\infty) \times Y )$, for $i=1,2$. Here we are considering that
the $T$ is related with the coordinate $u_i$ and the interval
$[0,\infty)$ with the coordinate $u_j$ for $i,j \in \{1,2\}$,
$i\neq j$ (see remark~\ref{rem: variables of cyl} below). All
these submanifolds are isometric in the Riemannian sense to $Z_i$
and we  identify their disjoint union with the Cartesian product
$Z_i \times [0,\infty)$

Let $E$ be a Hermitian vector bundle over a c.m.w.c.2, $X$. Let
$\Delta$ be a generalized Laplacian acting on $C^\infty(X,E)$, the
  sections of the vector bundle $E$.
The operator $\Delta$ is a {\bf compatible
  Laplacian} over $X$ if the following properties are satisfied:
\begin{itemize}
\item[i)] There exists a Hermitian vector bundle $E_i$ over $Z_i$
  such that $E \vert_{\R_+\times Z_i}$ is the pullback of $E_i$ under
  the projection $\pi:\R_+\times Z_i \to Z_i$, for $i=1,2$. We suppose
  also that the Hermitian metric of $E$ is the pullback of the
  Hermitian metric of $E_i$. On $\R_+ \times Z_i$, we have $
    \Delta=-\parcial{^2}{u_i^2}+\Delta_{Z_i}$,
  where $\Delta_{Z_i}$ is a compatible Laplacian acting on $C^\infty(Z_i,E_i)$. \item[ii)] There exists a Hermitian vector bundle $S$ over $Y$
such
  that $E \vert_{\R_+^2 \times Y}$ is the pullback of $S$ under the
  projection $\pi:\R_+^2 \times Y \to Y$. We assume also that the Hermitian product on $E \vert_{\R_+^2 \times Y}$ is the pullback of the Hermitian product on $S$. Finally we suppose that the operator $\Delta$ restricted to $\R_+^2 \times Y$  satisfies
  $
    \Delta=-\parcial{^2}{u_1^2}-\parcial{^2}{u_2^2}+\Delta_Y,
  $
  where $\Delta_Y$ is a generalized Laplacian acting on
  $C^\infty(Y,S)$.
\end{itemize}
Examples of compatible Laplacians  are given by the Laplacian
acting on forms and Laplacians associated to compatible Dirac
operators (see
 \cite{MuellerCorner}), they satisfy conditions i) and ii)
 due to the product structure of the Riemannian metric
 on the submanifolds $Y \times \R_+^2$ and $Z_i \times \R_+$. Since $X$ is a manifold with bounded
geometry and the vector bundle $E$ has bounded Hermitian metric,
the operator $\Delta:C^\infty_c(X,E)\subset L^2(X,E) \to L^2(X,E)$
is essentially self-adjoint (see~\cite[Corollary 4.2]{Shubin}).
Similarly $\Delta_{Z_i}:C^\infty_c(Z_i,E_i) \subset L^2(Z_i,E_i)
\to L^2(Z_i,E_i)$ is also essentially self-adjoint for $i=1,2$.
\begin{rem}\label{rem: variables of cyl}
If $j,k \in \{1,2\}$ and $j \neq k$, then we will denote by $u_j$
 the coordinate in $\R_+$ in the cylinder $Y \times \R_+$
of the complete manifold with cylindrical end $Z_k$.
\end{rem}
\begin{defin}\label{defin:hi Hi}
\begin{itemize}
\item Let $H$ and $H^{(i)}$ be the  self-adjoint extensions of
$\Delta:C^\infty_c(X,E)  \to L^2(X,E)$ and
$\Delta_{Z_i}:C^\infty_c(Z_i,E_i)  \to L^2(Z_i,E_i)$ respectively.
\item Let $b_i$ be the self-adjoint extension of
$-\frac{d^2}{du_i^2}:C^\infty_c(\R_+) \to L^2(\R_+)$ obtained by
imposing Dirichlet boundary conditions at $0$. \item Let $H_i$ be
the self-adjoint operator $b_i\otimes Id+Id \otimes H^{(i)}$
acting on $L^2(\R_+)\otimes L^2(Z_i,E_i)$. \item Let $H^{(3)}$ be
the self-adjoint operator associated to the essentially
self-adjoint operator $\Delta_Y:C^\infty(Y,S)\subset  L^2(Y,S)\to
L^2(Y,S)$ and let $H_3$ be the self-adjoint operator
$H_3:=b_1\otimes Id\otimes Id+Id\otimes b_2\otimes Id+ Id\otimes
Id\otimes H^{(3)}$ acting on $L^2(\R_+) \otimes L^2(\R_+) \otimes
L^2(Y,S)$. \item The operators $H_i$ are called {\bf channel
operators} for $i=1,2,3$.
    \end{itemize}
\end{defin}

The self-adjoint operators $H_1$ and $H_2$ have a free channel of
dimension 1 (associated to $b_1$ and $b_2$, respectively); the
operator $H_3$ has a free channel of dimension 2 (associated to
$b_1\otimes Id\otimes Id+Id\otimes b_2\otimes Id$).  In some parts
of this text we make an abuse of  notation by denoting $H$, $H_i$,
and $H^{(i)}$
 the
Laplacians acting on distributions and the self-adjoint operators
previously defined.

It is known that the compatible Laplacian $H^{(k)}$ decomposes the
Hilbert space $L^2(Z_k,E)$ into the orthogonal $H^{(k)}$-invariant
subspaces $L^2_{pp}(Z_k,E)$ and $L^2_{ac}(Z_k,E)$  associated to
pure point states  and absolutely continuous states (see
\cite{GUILLOPE} \cite{HUS}). We have $H^{(k)}=H^{(k)}_{pp}\oplus
H^{(k)}_{ac}$ on $L^2(Z_k,E)=L^2_{pp}(Z_k,E)\oplus
L^2_{ac}(Z_k,E)$ where $H^{(k)}_{pp}$ and $H^{(k)}_{ac}$ are
self-adjoint operators acting on  $L^2_{pp}(Z_k,E)$ and
$L^2_{ac}(Z_k,E)$. We define the self-adjoint operators
$H_{k,pp}:=b_k \otimes 1+1\otimes H^{(k)}_{pp}$ acting on
$L^2(\R_+) \otimes L^2_{pp}(Z_k,E_k)$, for $k=1,2$, that together
with $H$ will  define    important wave-operators in this article.
We notice that the operators $H_{k,pp}$ and $H^{(k)}_{pp}$ are
different operators, to see that we observe that they act in
different Hilbert spaces,  $H_{k,pp}$ has only absolutely
continuous spectrum, and $H^{(k)}_{pp}$ has only pure point
spectrum. Similarly, we define the self-adjoint operators $
H_{k,ac}:= b_k \otimes 1+1\otimes H^{(k)}_{ac} $ acting on
$L^2(\R_+) \otimes L^2_{ac}(Z_k,E_k)$. The operators $H_{k,ac}$
together with $H$ define   important wave--operators (see
theorem~\ref{thm:exist WO}).
\subsection{Main results}
Our first result is:
\begin{thm}\label{thm:exist WO}
1) For $k=1,2$ the following strong limits exist
$$W_{\pm}(H,H_{k,pp}):=\lim_{t \to \mp \infty} e^{ itH}e^{-
itH_{k,pp}}$$
$$W_{\pm}(H,H_{k,ac}):=\lim_{t \to \mp \infty} e^{
itH}e^{- itH_{k,ac}},$$ $$W_{\pm}(H,H_{3}):=\lim_{t \to \mp
\infty} e^{ itH}e^{- itH_{3}},$$ $$ W_\pm(H_{k,ac},H_3):=\lim_{t
\to \mp \infty} e^{ itH_{k,ac}}e^{-itH_3}.$$
\\
\\
2) The images of the operators $W_{\pm}(H,H_{1,pp})$,
$W_{\pm}(H,H_{2,pp})$ and $W_{\pm}(H,H_3)$ are pairwise
orthogonal.
\end{thm}
We call the operators defined in part 1) of the theorem  {\bf wave
operators}.
\begin{defin}\label{defin:asymp completen}We say that the wave operators, $W_{\pm}(H,H_{1,pp})$,
$W_{\pm}(H,H_{2,pp})$ and $W_{\pm}(H,H_3)$, are {\bf
asymptotically complete} if  for all $\psi \in L^2_{ac}(X,E)$
there exists $\varphi_k \in L^2_{pp}(Z_k,E_k)\otimes L^2(\R_+)$,
for $k=1,2$, and $\varphi_3 \in L^2(Y,S)\otimes L^2(\R_+^2)$ such
that
\begin{equation} \label{eq: asym completeness}
\psi=W_{\pm}(H,H_3) \varphi_3+\sum_{k=1}^2
W_{\pm}(H,H_{k,pp})\varphi_k.
\end{equation}
\end{defin}
Our second result is:
\begin{thm}\label{thm:asympt completeness}
The wave operators  $W_{\pm}(H,H_{1,pp})$, $W_{\pm}(H,H_{2,pp})$
and $W_{\pm}(H,H_3)$ are  asymptotically complete.
\end{thm}
Section~\ref{quantumdynam}   provides the first relation between
the quantum dynamics of the compatible Laplacian and the geometry
of $X$; theorem~\ref{thm:exist WO} is proved  in
section~\ref{existandort} using stationary phase methods. We prove
theorem~\ref{thm:asympt completeness} in
section~\ref{asymclustering} based on the methods of \cite{Y1}.
In~\ref{app: stationary phase} we give a summary of the stationary
phase methods used in section~\ref{existandort}.
\subsection{Related literature}
The literature about quantum scattering theory on open manifolds
is large. For that reason we restrict our bibliography to some
recent articles on the subject, where the reader can find
references to classic or basic articles, or to articles that we
consider directly related to the topics of this article. Articles
on quantum scattering theory on manifolds with cylindrical ends
are~\cite{ISOSCATTER},\cite{MUELLERSTRO}
 and \cite{RICHARDTIEDRA}; on manifolds asymptotically Euclidean
~\cite{MELeuclidean}; on $SL(3)/SO(3)$ \cite{MV3}; on homogeneous
spaces associated to finite groups on~\cite{BunkeOlbrich};
connections between scattering theory on compact asymptotically
Einstein manifolds and conformal geometry are studied
in~\cite{GrahamZworski}; quantum scattering theory on more general
open manifolds can be found in~\cite{CARRON} and~\cite{MUELSALOM}.
Relations between the geometry of manifolds with corners and the
quantum dynamics of many--particle Schr\"odinger operators has
been treated also in~\cite{VASY} but the topics are  different to
ours, and in particular the operators studied there are
many--particle Schr\"odinger operators that are essentially
perturbations via potentials of the Laplacian on $\R^n$, here we
treat perturbations associated to the geometry and not to a
potential. In~\cite{MuellerCorner} the spectral theory of
compatible Laplacians on c.m.w.c.2 is studied near $0$ under the
hypothesis that the compatible Laplacian on the corner has kernel
$0$, in this article we eliminate this hypothesis and study the
whole spectrum of the compatible Laplacians.
\section*{Aknowledgments}
This article is based on results obtained during the doctoral
studies of the author, realized under the supervision of Werner
M\"uller, at the University of Bonn and as well as work carried
out at Universidad de los Andes during the author's postdoc. He
wants to thank both institutions for the nice environments they
provided and Professor M\"uller  for his continuous support.
\section{Ruelle's theorem}\label{quantumdynam}
In this section we formulate Ruelle's theorem in the context of
compatible Laplacians on complete manifolds with a corner of
codimension 2, our aim is to give a first relation between the
quantum dynamics of the  compatible Laplacian and the geometry of
the manifold $X$.

Let $A$ be a self-adjoint operator acting on a Hilbert space
$\mathscr{H}$. We denote  $\mathscr{H}_{pp}(A)$  the subspace
spanned by all eigenvectors of $A$,
$\mathscr{H}_{c}(A):=(\mathscr{H}_{pp})^\perp(A)$,
$\mathscr{H}_{ac}(A)$, $\mathscr{H}_{sc}(A)$ will denote the
absolutely continuous and singular continuous subspaces of
$\mathscr{H}$ associated to $A$.
\begin{thm}{\upshape (cf. \cite[page 3452]{HS1} )}\label{thm: Ruelle}
Let $A$ be a self-adjoint operator acting  on $L^2(X,E)$ and
suppose that $A$ satisfies
\begin{equation}\label{eq: local prop}
\chi_K (A-\lambda)^{-1} \text{ is a compact operator, for any
compact subset $K$ of $X$,}
\end{equation}
for each $\chi_K \in C^\infty_c(X)$ such that $\chi_K=1$
restricted to $K$. Then:
\begin{equation*}
\begin{split}
&\varphi \in \mathscr{H}_{pp}(A) \Leftrightarrow \lim_{R \to \infty} \norm{(1-\chi_{R})e^{iAt}\varphi}=0 \text{ uniformly in $0 \leq t < \infty$.}\\
&\varphi \in \mathscr{H}_c(A) \Leftrightarrow \lim_{t \to \infty}
t^{-1}\int_0^t  \norm{\eta_{R} e^{iAs}\varphi}^2ds=0 \text{ for
any $R<\infty$,}
\end{split}
\end{equation*}
where  $\eta_R$ is any function in $C^\infty_c(X)$ that is equal
to $1$  on $X_R$, the compact manifold with a corner of
codimension 2 defined in~(\ref{eq: exhaustion X}).
\end{thm}
It follows from classical results in global analysis (see for
example~\cite{Shubin}) that the compatible Laplacian $H$
satisfies~(\ref{eq: local prop}). Then, intuitively, theorem
\ref{thm: Ruelle} implies that the  continuous states associated
to $H$ are moving away of compact sets as $t \to \infty$.
Theorems~\ref{thm:exist WO} and~\ref{thm:asympt completeness}
describe in more detail the asymptotic behavior of this escape.
\section{Existence of the wave operators}\label{existandort}
In this section we prove part 1) of theorem~\ref{thm:exist WO}
using Cook's criterion as expressed in the following simple lemma
of abstract scattering theory. We will make use also of stationary
phase methods which are summarized in~\ref{app: stationary phase}.
\begin{lem}\label{lem: existence WO abstract scattering yafaev}{ \upshape \cite[page 84]{Y3}}
 Let $B$ and $B_0$ be self-adjoint operators acting on Hilbert spaces
  $\mathscr{H}$ and $\mathscr{H}_0$ respectively. Let $\mathscr{J}:\mathscr{H}_0 \to \mathscr{H}$ be a bounded operator
 that takes the domain $\mathrm{Dom}(B_0)$ into the domain
 $\mathrm{Dom}(B)$. Suppose
 that for some $D_0 \subset \mathrm{Dom}(B_0)\cap \mathscr{H}_{0,ac}(B_0)$
 dense in
 $\mathscr{H}_{0,ac}(B_0)$, for any  $f \in D_0$
\begin{equation} \label{eq: existence wave operators}
\int_0^{\pm \infty} \norm{(B\mathscr{J}-\mathscr{J}B_0 )\exp(\mp
it B_0))f}dt < \infty.
\end{equation}
Then: $W_{\pm}(B,B_0,\mathscr{J}) :=s-\lim_{t \to \infty} \exp(\pm
i tB)\mathscr{J}\exp(\mp it B_0)$ exists.
\end{lem}
We prove first  the existence of $W_{\pm}(H,H_{k,pp})$, for $k \in
\{1,2\}$.

Let $\{\varphi_{k,j}\}_{j=1}^{N_k}$ be an orthonormal collection
of $L^2$--eigenfunctions of the operator $H^{(k)}_{pp}$ that
generates $L^2_{pp}(Z_k,E_k)$ for $k=1,2$. Observe that $N_1$ and
$N_2$ denote the number of $L^2$--eigenvalues of the Laplacians
$H^{(1)}$ and $H^{(2)}$ (counted with multiplicity). As pointed
out in~\cite{CANOTHESIS} and~\cite{CanoMourre},  the number of
$L^2$--eigenvalues of a Laplacian on a manifold with a cylindrical
end can be $0$, finite or infinite.  Without lost of generality
for our computations we will assume that there are infinite
$L^2$--eigenvalues that is $N_1=N_2=\infty$. Given $a \in
L^2(\R_+)$, $\hat{a}(u):=\int_0^\infty a(v)\sin v dv$ will denote
the sine transform of $a$. Let $\kappa \in C^\infty(\R_+)$ be such
that $\kappa(u)=0$ for $u\leq 2$ and $\kappa(u)=1$ for $u>3$. Let
us define $\kappa_k \in C^\infty(Z_k \times \R_+)$ by
$\kappa_k(z_k,u_k):=\kappa(u_k)$ for $k=1,2$ and extend it to
$C^\infty(X)$ by making it $0$ on $X-(Z_k \times \R_+)$. We will
show that we can apply lemma~\ref{lem: existence WO abstract
scattering yafaev} taking $\mathscr{J}=\kappa_k$, $B_0=H_{k,pp}$
and $B=H$. It is easy to see that $\kappa_k$ takes
$\mathrm{Dom}(H_{k,pp})$ into $\mathrm{Dom}(H)$. Let us denote by
$\mathscr{S}((0,\infty))$ the set of $C^\infty$--functions of
$[0,\infty)$ whose derivatives decrease faster than any polynomial
and such that all their derivatives at $0$ are equal to $0$. We
take
\begin{equation*}
\label{D0 for Hkpp} D_0:=\{g \varphi_{k,j}: j \in \N, g\in
\mathscr{S}((0,\infty)) \text{ and }\hat{g} \in
C^\infty_c((0,\infty))\}.\end{equation*} Since
$\mathscr{S}((0,\infty))$ is dense in $L^2(\R_+)$, it is easy to
see that the set $D_0$ is dense in $\mathrm{Dom}(H_{k,pp})$.
\\
\\
To prove~(\ref{eq: existence wave operators}) of lemma~\ref{lem:
existence WO abstract scattering yafaev} observe that for $f \in
D_0$
\begin{equation}\label{eq:ineq exist WHHkpp}
\begin{split}
&\norm{(H\kappa_k-\kappa_k H_{k,pp})e^{\mp itH_{k,pp}}f}\\
&\hspace{2cm} \leq \norm{\parcial{^2}{u^2_k}(\kappa_k)e^{\mp
itH_{k,pp}}f}+ 2
\norm{\parcial{}{u_k}(\kappa_k)\parcial{}{u_k}e^{\mp
itH_{k,pp}}f}. \end{split}
\end{equation}
If $f=g \varphi_{k,j} \in D_0$, we have
\begin{equation}\label{eq: bounding parcial eitHkkpp1}
\norm{\parcial{^2}{u^2_k}(\kappa_k)e^{\mp
itH_{k,pp}}f}=\norm{\frac{d^2}{du^2_k}(\kappa_k)e^{\mp
itb_k}g}_{L^2(\R_+)}.
\end{equation}
We can use~\ref{app: stationary phase} to see
$\int_{-\infty}^\infty \norm{\frac{d^2}{du^2_k}(\kappa_k)e^{\mp
itb_k}g}_{L^2(\R_+)}dt< \infty$. To estimate
$\parcial{}{u_k}(\kappa_k)\parcial{}{u_k}e^{\mp itH_{k,pp}}f$,
observe that
\begin{equation}\label{eq: bounding parcial eitHkkpp2}
\norm{\parcial{}{u_k}(\kappa_k)\parcial{}{u_k}e^{\mp
itH_{k,pp}}f}=\norm{\frac{d}{du_k}(\kappa_k)e^{\mp
itb_k}\frac{d}{du_k}g}_{L^2(\R_+)},
\end{equation}
then we can apply again the methods of~\ref{app: stationary
phase}. Finally, lemma~\ref{lem: existence WO abstract scattering
yafaev}  proves the existence of $W_{\pm}(H,H_{k,pp},\kappa_k)$.
\begin{prop}\label{prop:W(H,Hkpp) equal W(H,Hkpp,k)}
$W_{\pm}(H,H_{k,pp})$ exists and
$W_{\pm}(H,H_{k,pp},\kappa_k)=W_{\pm}(H,H_{k,pp})$.
\end{prop}
{\bf Proof:}
\\
Observe that for $f=g \varphi_{k,j} \in D_0$, we have $
\norm{e^{itH}(1-\kappa_k)e^{itH_{k,pp}}f}=\norm{(1-\kappa_k)e^{itb_k}g}_{L^2(\R_+,du_k)}.
$ Since $1-\kappa_k$ as a function of $u_k$ has compact support,
\ref{app: stationary phase}  implies $ s-\lim_{t \to
\infty}e^{itH}(1-\kappa_k)e^{itH_{k,pp}}=0.\square $

To prove the existence of $W(H, H_3)$ and $W(H, H_{k,ac})$ we
proceed analogously. Let $\{\phi_n\}_{n=0}^\infty$ be an
orthonormal collection of $L^2$--eigenfunctions of the operator
$H^{(3)}$ that generates $L^2(Y,S)$. We take  as dense sets
\begin{equation*}
\label{D0 for H3} D_{0,H_3}:=\{fg \phi_{n}: f\in
\mathscr{S}((0,\infty)_{u_1}), g\in \mathscr{S}((0,\infty)_{u_2})
\text{ and }\hat{f},\hat{g} \in
C^\infty_c((0,\infty))\},\end{equation*} and
\begin{equation*}
\label{D0 for Hac} D_{0,H_{k,ac}}:=\{f(z_k)g(u_k): f\in
\mathrm{Dom}(H^{(k)}_{ac}), g \in \mathscr{S}((0,\infty)) \text{
and }\hat{g} \in C^\infty_c((0,\infty))\}.
\end{equation*}
It is easy to see that~(\ref{eq:ineq exist WHHkpp})--(\ref{eq:
bounding parcial eitHkkpp2}) generalize and we can apply
lemma~\ref{lem: existence WO abstract scattering yafaev} to prove
the existence of $W(H, H_3, \kappa_1 \kappa_2)$ and $W(H,
H_{k,ac}, \kappa_k)$. Finally, there are natural  generalizations
of proposition~\ref{prop:W(H,Hkpp) equal W(H,Hkpp,k)} that show
the existence of $W(H, H_3)$ and $W(H, H_{k,ac})$.

The existence of $W_\pm(H_{1,ac},H_3)$ follows from the existence
of $W_\pm(H_{ac}^{(1)},b_2+H^{(3)})$  (see \cite{GUILLOPE}) and
the following equality
\begin{equation*}
W_{\pm}(b_1+H_{ac}^{(1)},b_1+b_2+H^{(3)})=Id_{L^2(\R_+,du_1)}
\otimes W_\pm(H_{ac}^{(1)},b_2+H^{(3)}).
\end{equation*}
\section{Orthogonality of the wave operators}\label{sec: ortho wave}
We prove part 2) of theorem~\ref{thm:exist WO}.
\subsection{Orthogonality of $W(H,H_{1,pp})$ and $W(H,H_{2,pp})$}
\label{sec:ortho waves Hkd}In this section we prove that for all
$f_k \in L^2_{pp}(Z_k,E_k) \otimes L^2(\R_+)$, $k=1,2$,
 the following equality holds
\begin{equation} \label{eq: orthogonality of Wave Zk}
\langle W_\pm(H, H_{1,pp})f_1,W_\pm(H,H_{2,pp})f_2
\rangle_{L^2(X,E)}=0.
\end{equation}
We observe that
\begin{equation*}
\begin{split}
& \langle W_\pm(H, H_{1,pp})f_1,W_\pm(H,H_{2,pp})f_2 \rangle_{L^2(X,E)}\\
&\hspace{2cm}=\lim_{t \to \infty} \langle e^{\mp i
tH_{1,pp}}f_1,e^{\mp itH_{2,pp}}f_2\rangle_{L^2(X,E)},
\end{split}
\end{equation*}
hence, equation~(\ref{eq: orthogonality of Wave Zk}) is satisfied
as a consequence of the following lemma.
\begin{lem}\label{lem:orthogon. using eiHkdt}
For all $f_k \in L^2_{pp}(Z_k,E_k) \otimes L^2(\R_+)$, $k=1,2$,
\begin{equation*} \label{eq: orthogon. using eiHkdt}
\lim_{t \to \infty} \langle e^{\mp i tH_{1,pp}}f_1,e^{\mp
itH_{2,pp}}f_2\rangle_{L^2(X,E)}=0.
\end{equation*}
\end{lem}
{\bf Proof:}
\\
By continuity of the bilinear form
$$(f_1,f_2) \mapsto \langle W_\pm(H, H_{1,pp})f_1,W_\pm(H,H_{2,pp})f_2\rangle_{L^2(X,E)},$$
it is enough to prove the lemma for the dense set of functions of
the form $f_k=a_k \varphi_k$, where $\varphi_k \in L^2(Z_k,E_k)$
is an $L^2$-eigenfunction of $H^{(k)}$ with eigenvalue $\gamma_k$,
$a_k \in \mathscr{S}((0,\infty))$ and $\hat{a}_k\in
C^\infty_c((0,\infty))$.

In the next computation we use the notation given in
definition~\ref{defin:hi Hi}  and explained in remark~\ref{rem:
variables of cyl},
\begin{equation}\label{eq:calc ortho WHkd 2}
\begin{split}
&\normv{\langle e^{\mp i tH_{1,pp}}f_1, e^{\mp itH_{2,pp}}f_2
\rangle_{L^2(X,E)}}\leq \int \vert \langle \int_0^\infty
\varphi_1(u_2,y) \cdot e^{\pm i tb_2}
(a_2)(u_2)du_2,\\
&\hspace{4cm}\int_0^\infty \varphi_2(u_1,y) \cdot e^{\pm i
tb_1}(a_1)(u_1) du_1 \rangle \vert dvol(y),
\end{split}
\end{equation}
where the Hermitian product inside the integrals on the
right--hand side of the inequality is the Hermitian product of the
vector bundle $S \to Y$.   It is well known that there exists $C
\in \R$ such that $\normv{e^{\pm i tb_k} (a_k)(u_k)} \leq C
t^{-1/2},$ for all $t>1$ and for all $u_k \in \R_+$
(see~\cite[Corollary, page 41]{RS3}). Cauchy--Schwartz applied to
the last term of~(\ref{eq:calc ortho WHkd 2}) and the fact
$\normv{\varphi_k(u_j,y)} \leq C e^{-cu_1}$
 for some $c>0$ (see \cite[Lemma 1.36]{HUS})
 finish the proof of the lemma.
$\square$
\subsection{$\mathrm{Im}(W_{\pm}(H,H_3))$ is orthogonal to $\mathrm{Im}(W_\pm(H,H_{k,pp}))$}
\label{sec:ortho omega y Hkd} Without lost of generality we prove
the orthogonality of $\mathrm{Im}(W_{\pm}(H,H_3))$ and
$\mathrm{Im}(W_\pm(H,H_{1,pp}))$. Let $\phi \in L^2(Y,S)$,
$\varphi \in L^2_{pp}(Z_1,E_1)$, $c\in L^2(\R_+,du_1)$ and $a_i
\in L^2(\R_+,du_i)$ for $i=1,2$. It is enough to prove that
\begin{equation}\label{eq:ortho H3 yHpp}
\langle W_\pm(H, H_3)(a_1a_2 \phi),W_\pm(H,H_{1,pp})(c\varphi)
\rangle_{L^2(X,E)}=0. \end{equation} We have
\begin{equation}\label{eq:ortho H3 y Hkpp 2}
\begin{split}
&\normv{ \langle e^{\pm itH_3}(a_1a_2 \phi),e^{\pm
itH_{1,pp}}(c\varphi) \rangle_{L^2(X,E)}}\\
&\hspace{1cm}\leq \normv{\langle e^{\pm
it(b_2+H^{(3)})}(a_2\phi),e^{\pm itH^{(1)}_{pp}}(\varphi)
\rangle_{L^2(Z_1,E_1)}}.
\end{split}
\end{equation}
Since the wave operator $W_\pm(H^{(1)},b_2+H^{(3)})$ is complete,
we can find $\psi \in L^2_{ac}(Z_1,E_1)$ such that
$$\lim_{t \to \infty}\norm{e^{\pm it(b_2+H^{(3)})}(a_2\phi)-e^{\pm itH^{(1)}}\psi
}_{L^2(Z_1,E_1)}=0.$$ This together with~(\ref{eq:ortho H3 y Hkpp
2}) imply~(\ref{eq:ortho H3 yHpp}), since $ L^2_{ac}(Z_1,E_1)$ is
orthogonal to $ L^2_{pp}(Z_1,E_1)$.
\section{Asymptotic clustering:
a time dependent approach} \label{asymclustering}In this section
we prove asymptotic completeness (theorem~\ref{thm:asympt
completeness}) using a time dependent approach. We follow
closely~\cite{HS1}, and as in this article our main tools will be
Mourre's inequality and the Yafaev functions (see section
\ref{section: Yafaev function}) properly adapted to the context of
compatible Laplacians on c.m.w.c.2.
\subsection{Mourre estimate for compatible Laplacians}\label{Mourre estimate}
First we state  Mourre's inequality which will be used to prove
asymptotic completeness. It was developed in~\cite{CanoMourre} and
used
 to prove the absence of singular continuous
spectrum of compatible Laplacians on c.m.w.c.2 and also to prove
that the pure point spectrum of these operators accumulates only
at thersholds. Let $\kappa \in C^\infty(\R_+)$ be such that
$\kappa(u)=0$ for $u\leq 2$ and $\kappa(u)=1$ for $u>3$. Let us
define $\kappa_k \in C^\infty(Z_k \times \R_+)$ by
$\kappa_k(z_k,u_k):=\kappa(u_k)$ for $k=1,2$ and the function $r^2
\in C^\infty(\R_+^2)$ by $r^2(u_1,u_2):=\kappa(u_1) u_1^2
+\kappa(u_2) u_2^2$. The function $r^2$ induces a function on $Y
\times \R_+^2$ by $(y,u_1,u_2) \mapsto r^2(u_1,u_2) $ and this
function extends naturally to $X$ by making it $0$ out of $Y\times
\R_+^2$, by an abuse of notation we denote this new function by
$r^2$ too. We extend $\kappa_1$ and $\kappa_2$ to $X$ similarly by
making them $0$ out of $Z_1 \times \R_+$ and $Z_2 \times \R_+$
respectively. Let us define the first order differential operator
$A$  by
\begin{equation*}\label{eq:def A corner}
A:=i[H,r^2].
\end{equation*}
We define {\bf the set of thresholds of $H$}, $\tau(H)$, by
\begin{equation*}\label{def: thresholds}
\begin{split}
\tau(H):=\sigma_{pp}\left(H^{(1)}\right)\cup
\sigma_{pp}\left(H^{(2)}\right)\cup
\sigma_{pp}\left(H^{(3)}\right).
\end{split}
\end{equation*}
Let $\Sigma:=\min \tau(H)$, such a minimum exists  because
$H^{(1)}, H^{(2)}$ and $H^{(3)}$ are bounded from below
(see~\cite[Satz 1.27]{HUS}) and hence the three sets on the right
are discrete and with a minimum. For $\lambda \in \R$, define the
number
\begin{equation*}\label{eq:def theta}
\begin{split}
&\theta(\lambda):=
\begin{cases}
0,&\text{ for } \lambda \leq \Sigma;\\
\inf\{\lambda-\gamma:\gamma \in \tau(H), \gamma< \lambda\} ,&
\text{ for } \lambda > \Sigma.
\end{cases}
\end{split}
\end{equation*}
The next theorem is the generalization of Mourre's inequality to
c.m.w.c.2 that was developed in \cite{CanoMourre}.
\begin{thm}{\upshape \cite[theorem 5]{CanoMourre}}\label{thm:Mourre estimate}
Given $\lambda \in \R$ and $\varepsilon>0$, there exist an open
interval $I \ni\lambda$,  and  an $H$-compact operator $K$ such
that
  \begin{equation*}
    E_I(H)\,i[H,A]\,E_I(H)\geq (\theta(\lambda)-\varepsilon)\,E_I(H)+K,
  \end{equation*}
  where $E_I(H)$ denotes the spectral projection of the operator $H$ on the interval $I \subset \R$.
\end{thm}
\subsection{Graf-Yafaev functions}
\label{section: Yafaev function}Consider the Schr\"odinger
operators $\sum_{i=1}^2 \left(-\parcial{^2}{u^2_i} +V_i\right)$
acting on $L^2(\R^2)$ where $V_i \in C^\infty(\R^2)$, $V_i$
depends only of the variable $u_i$ and is compactly supported in
this variable. Our Graf-Yafaev functions  are constructed in
analogy to the Graf--Yafaev functions associated to these
Schr\"odinger operators following~\cite{HS1},~\cite{HS2}
and~\cite{Y1}. In this section we will omit some  proofs, because
we consider that the analogy is direct once the Graf--Yafaev
functions are constructed.

Given $\epsilon>0$, we take $\varepsilon^-_0:= l_0 \leq
\varepsilon_0 \leq l_0+\epsilon=:\varepsilon^-_0$,
$\varepsilon^-_3:=2 \epsilon^2 < \varepsilon_3 <
3\epsilon^2=:\varepsilon^+_3$, and $\varepsilon^-_i:=2 \epsilon <
\varepsilon_i < 3\epsilon:=\varepsilon^+_i$ for $i=1,2$. We call
the vectors $
\varepsilon:=(\varepsilon_1,\varepsilon_2,\varepsilon_3)$ {\bf
$\epsilon$--admissible}. From now on we will denote
$\normv{(u_1,u_2)}:=\sqrt{u^2_1+u^2_2}$.
\\
\\
Let $\chi$ be the characteristic function of the interval
$[0,\infty)$. The functions $g^{(1)},g^{(2)}$ and $g^{(3)}$
defined below are analogous to the functions $m^{(a)}$ in
\cite[equation 3.9]{Y1}.
\begin{equation*}
\begin{split}
&g^{(0)}({\bf \varepsilon},x):=\\
&\begin{cases}
\varepsilon_0 \chi\left(\varepsilon_0 -\max \{(1+\varepsilon_1)u_1,(1+\varepsilon_2)u_2,(1+\varepsilon_3)\normv{u}\}\right) \\
\hspace{1cm}\text{ for } x=(y,u_1,u_2) \in Y \times \R_+^2.\\
\varepsilon_0\hspace{1cm} \text{ if } x=(u_1,z_1) \in
[0,\frac{\varepsilon_0}{1+\varepsilon_1}] \times Z_{1,0}.\\
\varepsilon_0\hspace{1cm} \text{ if } x=(u_2,z_2) \in
[0,\frac{\varepsilon_0}{1+\varepsilon_2}] \times Z_{2,0}.\\
\varepsilon_0\hspace{1cm} \text{ if } x\in X_0\\
0 \hspace{1cm} \text{otherwise}.
\end{cases}
\end{split}
\end{equation*}
\begin{equation*}
\begin{split}
&g^{(1)}({\bf \varepsilon},x):=\\
&\begin{cases}
(1+\varepsilon_1)u_1 \chi\left((1+\varepsilon_1)u_1 -\max \{\varepsilon_0,(1+\varepsilon_2)u_2,(1+\varepsilon_3)\normv{(u_1,u_2)}\}\right) \\
\hspace{1.5cm}\text{ for } x=(y,u_1,u_2) \in Y \times \R_+^2.\\
(1+\varepsilon_1)u_1 \hspace{0.8cm} \text{ for } x=(z_1,u_1) \in Z_{1,0} \times [\frac{\varepsilon_0}{1+\varepsilon_1},\infty).\\
0 \hspace{1.3cm}\text{ otherwise}.
\end{cases}
\end{split}
\end{equation*}
\begin{equation*}
\begin{split}
&g^{(2)}({\bf \varepsilon},x):=\\
&\begin{cases}
(1+\varepsilon_2)u_2 \chi\left((1+\varepsilon_2)u_2 -\max \{\varepsilon_0,(1+\varepsilon_1)u_1,(1+\varepsilon_3)\normv{(u_1,u_2)}\}\right) \\
\hspace{1.5cm}\text{ for } x=(y,u_1,u_2) \in Y \times \R_+^2.\\
(1+\varepsilon_2)u_2 \hspace{0.8cm} \text{ for } x=(z_2,u_2) \in Z_{2,0} \times [\frac{\varepsilon_0}{1+\varepsilon_2},\infty).\\
0 \hspace{1.3cm} \text{ otherwise}.
\end{cases}
\end{split}
\end{equation*}
\begin{equation*}
\begin{split}
&g^{(3)}({\bf \varepsilon},x):=\\
&\begin{cases}
(1+\varepsilon_3)\normv{(u_1,u_2)} \chi\left((1+\varepsilon_3)\normv{(u_1,u_2)} -\max \{\varepsilon_0,(1+\varepsilon_2)u_2,(1+\varepsilon_1)u_1\}\right) \\
\hspace{1.5cm}\text{ for } x=(y,u_1,u_2) \in Y \times \R_+^2.\\
0 \hspace{1.3cm}\text{ otherwise}.
\end{cases}
\end{split}
\end{equation*}
The functions $g^{(i)}(x,\varepsilon)$ could be defined more
directly in our case, for example for $(y,u_1,u_2) \in Y \times
\R_+^2$, $\varepsilon$ $\epsilon$--admissible and $\epsilon$ small
enough, $g^{(1)}({\bf
\varepsilon},y,u_1,u_2)=(1+\varepsilon_1)u_1$ if and only if
$(1+\varepsilon_1)u_1$ is greater or equal than $\varepsilon_0$,
$(1+\varepsilon_2)u_2$ and $(1+\varepsilon_3)\normv{(u_1,u_2)}$,
and $g^{(1)}({\bf \varepsilon},y,u_1,u_2)=0$ otherwise. However we
used the previous definitions because we would like to point out
that the Graf--Yafaev method could generalize to manifolds with
corners of higher codimension. Let us define the function
\begin{equation*}
\begin{split}
&g(x,\varepsilon):=
\begin{cases}
(1+\varepsilon_i)u_i \hspace{0.5cm}\text{ for } x=(z_i,u_i) \in Z_{i,0} \times \R_+.\\
\max\{\varepsilon_0,(1+\varepsilon_1)u_1,(1+\varepsilon_2)u_2,(1+\varepsilon_3)\normv{u}\},\\
\hspace{0.8cm}\text{ for } x=(y,u_1,u_2) \in Y \times \R_+^2.\\
\varepsilon_0 \hspace{1cm} \text{ for } x \in X_0.
\end{cases}
\end{split}
\end{equation*}
We observe that
\begin{equation}\label{eq:g in terms gi}
g(x,\varepsilon)=\sum_{i=0}^3 g^{(i)}(x,\varepsilon).
\end{equation}
The next functions will be important in the description of the
functions $g$ and $g^{(i)}$.
\begin{equation*}
\begin{split}
&k_1(\varepsilon_1,\varepsilon_3):=\frac{1+\varepsilon_3}{\sqrt{(1+\varepsilon_1^2)-(1+\varepsilon_3)^2}},
\hspace{0.5 cm}
k_2(\varepsilon_1,\varepsilon_2):=\frac{1+\varepsilon_2}{1+\varepsilon_1}\\
&\hspace{2cm}\text{ and }
k_3(\varepsilon_2,\varepsilon_3):=\frac{\sqrt{(1+\varepsilon_2)^2-(1+\varepsilon_3)^2}}{1+\varepsilon_3}.
\end{split}
\end{equation*}
The next proposition   is a consequence of the following limits $
\lim_{\epsilon \to 0} k_1(\varepsilon_1,\varepsilon_3)=\infty$,
$\lim_{\epsilon \to 0} k_2(\varepsilon_1,\varepsilon_2)=1$, and
$\lim_{\epsilon \to 0} k_3(\varepsilon_2,\varepsilon_3)=0.
$\begin{prop}\label{Prop:lema defining gi} Let $\epsilon>0$ be
small enough and let
$\varepsilon:=(\varepsilon_1,\varepsilon_2,\varepsilon_3)$ be an
$\epsilon$--admissible vector. Then
\begin{equation*}
k_1(\varepsilon_1,\varepsilon_3) \geq
 k_2(\varepsilon_1,\varepsilon_2) \geq
k_3(\varepsilon_2,\varepsilon_3).
\end{equation*}
\end{prop}
Proposition~\ref{Prop:lema defining gi} implies that
$(1+\varepsilon_1)u_1 \geq \max\{\varepsilon_0,
(1+\varepsilon_2)u_2,(1+\varepsilon_3)\normv{(u_1,u_2)} \}$ if and
only if $u_1 \geq \frac{\varepsilon_0}{1+\varepsilon_1}$ and $u_1
> k_1(\varepsilon_1,\varepsilon_3)u_2$. Reasoning in this way
 we
obtain the sketch of the function $g(x,\varepsilon)$ given in
figure 3.
\begin{center}
      \psset{unit=0.5cm}
          \begin{pspicture}(-3,-5)(10,10)
                  \psline(0,0)(0,10) 
                   \qline(0,0)(10,0)
        \psline(4,0)(4,3)
        \rput(-2,3){$u_1=\frac{\varepsilon_0}{1+\varepsilon_1}$}
        \psline(0,3)(2,3)
        \rput{45}(2,2.2){$g^{(0)}(\varepsilon,x)=\varepsilon_0$}
        \psarc[](2,2){3.2}{25}{75}
        \psline[linewidth=1.5pt,linestyle=dotted,dotsep=2pt](4,3)(10,6)
        \psline[linewidth=1.5pt,linestyle=dashed](2,3)(4.5,10)
        \psline[linewidth=1.5pt,linestyle=dashed](-3,-1)(-1.5,-1)
        \psline[linewidth=1.5pt,linestyle=dotted](-3,-2.5)(-1.5,-2.5)
        \rput[l](1,-1){$u_1 = k_1(\varepsilon_1,\varepsilon_3)u_2$}
        \rput[l](1,-2.5){$u_1 =k_3(\varepsilon_2,\varepsilon_3)u_2$}
        \rput{68}(2,7){$(1+\varepsilon_1) u_1=g^{(1)}(\varepsilon,x)$}
        \rput{35}(7,7){$(1+\varepsilon_3)\normv{u}=g^{(3)}(\varepsilon,x)$}
        \rput(8,2){$(1+\varepsilon_2) u_2=g^{(2)}(\varepsilon,x)$}
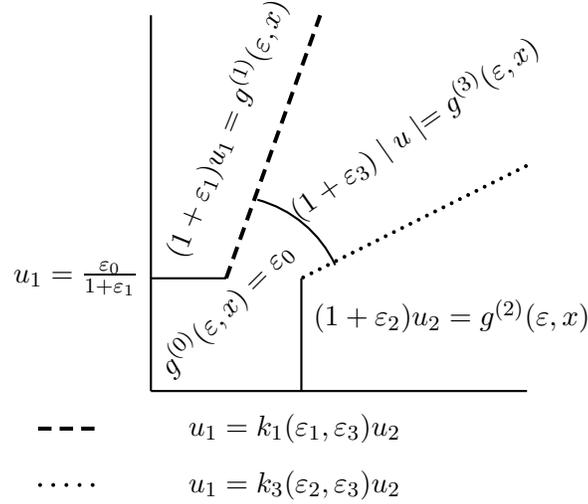
        \rput(3.5,-4){Figure 3. Sketch of the Graf--Yafaev function $g(x,\varepsilon)$.}
    \end{pspicture}
\end{center}
Let $\varphi_i \geq 0$, $\varphi_i \in C^\infty_c(\R_+)$,
$\mathrm{supp} \varphi_i \subset
[\varepsilon_i^-,\varepsilon_i^+]$ and $\int_0^\infty
\varphi_i(\varepsilon_i)d\varepsilon_i=1$, for $i=1,2,3$.  Let
$\varphi_0 \in C^\infty(\R_+)$ be   a real function with support
in the interval $(l_0,l_0 +\epsilon,)$ for some $l_0>0$, that
satisfies also $\int_0^\infty
\varphi_0(\varepsilon_0)d\varepsilon_0=1$. We regularize the
function $g^{(i)}$ averaging over the $\epsilon$--compatible
vectors $\varepsilon$:
\begin{equation}\label{eq:regularize gi}
\begin{split}
g^{(i)}(x):=\int_{-\infty}^\infty g^{(i)}(x, {\bf
\varepsilon})\Pi_{i=0}^3\left(\varphi_i(\varepsilon_i)d\varepsilon_i\right).
\end{split}
\end{equation}
Definition (\ref{eq:regularize gi}) is inspired by
\cite[definition 3.12]{Y1}. For $i=0,1,2,3$, define $
\Phi_i(\xi):= \int_0^\xi \varphi_i(\varepsilon_i)d\varepsilon_i$.
An easy computation shows
\begin{equation}\label{eq: def g1 int of epsilon}
\begin{split}
&g^{(1)}(x)=\int_{-\infty}^\infty (1+\varepsilon_1)u_1 \varphi_1(\varepsilon_1) \Phi_0\left((1+\varepsilon_1)u_1\right)\Phi_2\left((1+\varepsilon_1)u_1 u_2^{-1}-1\right)\\
&\hspace{1.5cm} \cdot
\Phi_3\left((1+\varepsilon_1)u_1\normv{(u_1,u_2)}^{-1}-1\right)
d\varepsilon_1,
\end{split}
\end{equation}
for $x=(y,u_1,u_2) \in Y \times \R_+^2$ or  $x=(z_1,u_1)\in Z_1
\times \R_+$. We observe that $g^{(1)}(x)=0$ on $X-(Z_1 \times
\R_+)$. There is a similar formula for $g^{(2)}(x)$. For $g^{(3)}$
and $g^{(0)}$ we have:
\begin{equation}\label{eq: def g3 int of epsilon}
\begin{split}
&g^{(3)}(x)=\int (1+\varepsilon_3)\normv{u} \varphi_3(\varepsilon_3) \Phi_0\left( (1+\varepsilon_3)\normv{u}\right)\Phi_1\left(\frac{(1+\varepsilon_3)\normv{u}}{ u_1}-1\right) \\
&\hspace{1.5cm}
\cdot\Phi_2\left(\frac{(1+\varepsilon_3)\normv{u}}{
u_2}-1\right)d\varepsilon_3,
\end{split}
\end{equation}
\begin{equation}\label{eq: def g0 int of epsilon}
\begin{split}
&g^{(0)}(x)=\int \varepsilon_0 \varphi_0(\varepsilon_0)
\Phi_1(\frac{\varepsilon_0}{u_1}-1)\Phi_2(\frac{\varepsilon_0}{u_2}-1)\Phi_3(\frac{\varepsilon_0}{\normv{u}}-1)d\varepsilon_0,
\end{split}
\end{equation}
for $x=(y,u_1,u_2) \in Y \times \R_+^2$.

We define $g$, the regularization of the function
$g(x,\varepsilon)$, by taking the average on $\varepsilon$ of
$g(x,\varepsilon)$.
\begin{equation*}
\begin{split}
g(x):=&\int \max\{\varepsilon_0,(1+\varepsilon_1)u_1,(1+\varepsilon_2)u_2, (1+\varepsilon_3)\normv{u}\}\\
&\varphi_0(\varepsilon_0)\varphi_1(\varepsilon_1)\varphi_2(\varepsilon_2)\varphi_3(\varepsilon_3)d\varepsilon_0d\varepsilon_1
d\varepsilon_2 d\varepsilon_3.
\end{split}
\end{equation*}
Let us define $\mu_0:=\int \varepsilon_0 \varphi_0(\varepsilon_0)
d\varepsilon_0$ and $\mu_i:=\int (1+\varepsilon_i)
\varphi_i(\varepsilon_i) d\varepsilon_i$ for $i=1,2,3$. We observe
that the maximum of the function $k_1$ in
$[2\epsilon,3\epsilon]\times [2\epsilon^2,3\epsilon^2]$ is
attained in
$(\varepsilon_1,\varepsilon_2)=(2\epsilon,3\epsilon^2)$ and its
minimum is attained in
$(\varepsilon_1,\varepsilon_2)=(3\epsilon,2\epsilon^2)$. The
maximum of the function $k_2$ in $[2\epsilon,3\epsilon]\times
[2\epsilon^2,3\epsilon^2]$ is attained in
$(\varepsilon_2,\varepsilon_3)=(3\epsilon,2\epsilon^2)$ and its
minimum is obtained in
$(\varepsilon_1,\varepsilon_2)=(2\epsilon,3\epsilon)$.  Based on
these observations and proposition~\ref{Prop:lema defining gi} we
obtain figure 4, a sketch of the Yafaev function. The arcs in this
figure are part of the circles
$\normv{(u_1,u_2)}=\frac{l_0}{1+3\epsilon^2}$ and
$\normv{(u_1,u_2)}=\frac{l_0+\epsilon}{1+2\epsilon^2}$.
\begin{center}
      \psset{unit=0.5cm}
          \begin{pspicture}(-3,-5)(10,10)
              \psline(0,0)(0,10) 
                   \qline(0,0)(10,0)
                \psline(0,3)(2.5,3)
                \rput(-2,3){$u_1=\frac{l_0}{1+3\epsilon}$}
                \psline(3,0)(3,2.5)
                \rput{90}(3,-2){$u_2=\frac{l_0}{1+3\epsilon}$}
                \psline(0,4)(3,4)
                \rput(-2,4){$u_1=\frac{l_0+\epsilon}{1+2\epsilon}$}
                \psline(4,0)(4,3)
                \rput{90}(4.5,-2){$u_2=\frac{l_0+\epsilon}{1+2\epsilon}$}
                \rput{48}(2.1,2.1){$g^{(0)}(x)=\mu_0$}
                \rput{75}(1.5,7.7){$\mu_1 u_1=g^{(1)}(x)$}
        \rput{40}(7.5,7.5){$\mu_3\normv{u}=g^{(3)}(x)$}
        \rput(8,1){$\mu_2 u_2=g^{(2)}(x)$}
        \pspolygon[fillstyle=solid,fillcolor=lightgray](0,3)(2.2,3)(2.5,4)(0,4)
        \pspolygon[fillstyle=solid,fillcolor=lightgray](2.2,3)(2.5,3)(3,4)(2.5,4)(2.2,3)
        \pspolygon[fillstyle=solid,fillcolor=lightgray](2.5,4)(3,4)(6,10)(4,10)(2.5,4)
        \pspolygon[fillstyle=solid,fillcolor=lightgray](3,0)(3,2.2)(4,2.5)(4,0)
        \pspolygon[fillstyle=solid,fillcolor=lightgray](3,2.2)(3,2.5)(4,3)(4,2.5)(3,2.2)
        \pspolygon[fillstyle=solid,fillcolor=lightgray](4,2.5)(4,3)(10,6)(10,4)(4,2.5)
        \pscustom[fillstyle=solid,fillcolor=lightgray]{\psarc[](2,2){3.2}{26}{64}
            \psarcn[](2,2){3.8}{63}{27}\psarc[](2,2){3.2}{26}{64}}
            \psline[linewidth=2pt]{->}(8,5)(8,5.8)
            \pscircle(8,6.2){0.4}
            \rput(8,6.2){$3$}
            \psline[linewidth=2pt]{->}(5.5,9)(6.4,9)
            \pscircle(6.9,9){0.4}
            \rput(6.9,9){$2$}
            \psline[linewidth=2pt]{->}(3.4,7.5)(5.5,7.5)
            \pscircle(6,7.5){0.4}
            \rput(6,7.5){$1$}
            \psline[linewidth=2pt]{->}(9,3.7)(9,3)
            \pscircle(9,2.5){0.4}
            \rput(9,2.5){$4$}
            \rput(-4.5,-1){$1$}
            \pscircle(-4.5,-1){0.4}
            \rput(-1,-1){$u_1=k_1(2\epsilon,3\epsilon^2)u_2$}
            \rput(-4.5,-2){$2$}
            \pscircle(-4.5,-2){0.4}
            \rput(-1,-2){$u_1=k_1(3\epsilon,2\epsilon^2)u_2$}
            \rput(7,-1){$3$}
            \pscircle(7,-1){0.4}
            \rput(10.5,-1){$u_1=k_2 (3\epsilon,2\epsilon^2)u_2$}
            \rput(7,-2){$4$}
            \pscircle(7,-2){0.4}
            \rput(10.5,-2){$u_1=k_2 (2\epsilon,3\epsilon^2)$}
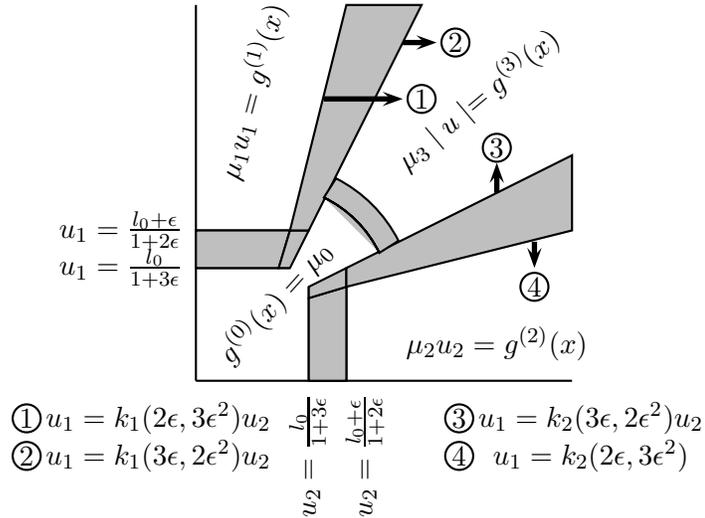
            \rput(3.5,-5){Figure 4. Sketch of the Graf--Yafaev functions.}
\end{pspicture}
\end{center}
The next lemma summarize the main properties of $g$ that we will
use in this article.
\begin{lem}{\upshape (cf.  \cite[page 538]{Y1})}\label{thm: properties of g}
$g$ satisfies the following properties:
\begin{itemize}
\item[1)] $g \in C^\infty(X)$ and $g(x)$ is real homogeneous of
degree $1$ in the sense that:
\begin{equation*}
\begin{split}
&g(tu_1,z_1)= tg(u_1,z_1) \text{ for } z_1 \in Z_{1}, u_1 \geq 4; \text{ and,}\\
&g(tu_1, tu_2, y)=t g(u_1, u_2, y) \text{ for } (y,u_1,u_2) \in Y
\times [4,\infty)^2
\end{split}
\end{equation*}
for $t \geq 0$ \item[2)] $g(x)\geq 1$ for $x \in X-X_{4}$.
\item[3)] $g(x)$ is convex in the sense that for $(y,u_1,u_2)$ and
$(y,v_1,v_2)$ $Y \times \R_+^2$:
\begin{equation*}
g((y,t(u_1,u_2)+s(v_1,v_2))\leq sg((y,u_1,u_2))+t g(y,v_1,v_2),
\end{equation*}
for  $s,t \in [0,1], s+t=1$. \item[4)] The functions $g^{(i)}$'s
are related to $g$ by the equality $g(x)=\sum_{i=0}^3 g^{(i)}(x)$.
\end{itemize}
\end{lem}
{\bf Proof:}
\\
(\ref{eq: def g1 int of epsilon}), (\ref{eq: def g3 int of
epsilon}) and (\ref{eq: def g0 int of epsilon}) prove that the
functions $g^{(i)}$ are smooth. 4) follows from (\ref{eq:g in
terms gi}) and these results imply that $g$ is smooth. $g$ is
convex because it is the integral of the maximum of convex
functions. The other properties follow from direct calculations.
$\square$
\begin{defin}
A function $g$ satisfying properties 1),2), 4) and 5) of the above
lemma is called {\bf Yafaev function}.
\end{defin}
Let $f:X \to \R$ be a $C^\infty$- function. Let us denote by $f''$
the matrix valued function
\begin{equation}\label{eq:defin Hessian}
f'':=\begin{pmatrix}
\parcial{^2}{u_1^2}f&\parcial{^2}{u_1\partial u_2}f\\
\parcial{^2}{u_1\partial u_2}f&\parcial{^2}{u_2^2}f
\end{pmatrix},\end{equation}
defined on $Y \times \R_+^2$. We observe that the matrix of
functions $(g^{(i)})''$ can be extended to $X$  making it $0$ out
of $Y\times \R^2_+$; we will make this type of natural extension
without to explicitly mention them for other functions. We remark
that $(f)''$ is not the Hessian of $f$.

According to the previous lemma, the functions $g^{(i)}$ are
Yafaev functions, but they do not satisfy 3), yet in any case they
are bounded by  suitable convex functions, as it is shown in the
next lemma.
\begin{lem}{\upshape (cf. \cite[lemma 7.4]{HS1})}\label{lem: bounding gi by a convex Yafaev}
For each $i \in \{1,2,3\}$ there exists $\tilde{g}_i$ a Yafaev
function  such that $ (\tilde{g}^{(i)})''(x) \geq (g^{(i)})''(x)$,
for all $x \in X$.
\end{lem}
{\bf Proof:} We prove the lemma for $i=1$, the  cases $i=2$  and
$i=3$ can be  treated similarly. Let us define the set
\begin{equation*}
\begin{split}
\Gamma:=&\{(y,u_1,u_2) \in Y \times \R_+^2:
\frac{l_0}{1+2\epsilon}
\leq u_1 \leq \frac{l_0+\epsilon}{1+3\epsilon} \text{ and }\\
&\hspace{1cm} k_1(3\epsilon,2\epsilon^2)u_2 \leq u_1 \leq
k_1(2\epsilon,3\epsilon^2)u_2\}.\end{split}
\end{equation*}
Let $x_0 \in \Gamma$ and let $\delta_{x_0}$ be a positive function
in $C^\infty_c(\R)$ such that $\delta_{x_0}(g(x_0)) \neq 0$.
Taking $\epsilon$ small enough and $l_0$ suitable we can find a
Yafaev function $g$ such that $(g)''(x_0)=(r)''(x_0)>0$. Let us
define the function
\begin{equation*}
\begin{split}
\tilde{g}_{x_0}(x)&:=\int_{g(x)}^{\infty}  s
\delta_{x_0}(s)ds+g(x)\int_{-\infty} ^{g(x)}\delta_{x_0}(s)ds.
\end{split}
\end{equation*}
We have
\begin{equation*}
\begin{split}
\tilde{g}''_{x_0}(x)=&g''(x)\int_{-\infty}^{g(x)}\delta_{x_0}(s)ds\\
&+\delta_{x_0}(g(x))\begin{pmatrix}
\parcial{}{u_1}(g)^2&\parcial{}{u_1}(g)\parcial{}{u_2}(g)\\
\parcial{}{u_1}(g)\parcial{}{u_2}(g)&\parcial{}{u_2}(g)^2
\end{pmatrix}(x).
\end{split}
\end{equation*}
Since $(g)''(x)$ is positive and
$\int_{-\infty}^{g(x)}\delta_{x_0}(s)ds>0$ for $x \in \Gamma$ near
enough to $x_0$, we have:
\begin{equation*}
\begin{split}
&\hspace{1cm}\langle\tilde{g}''_{x_0}(x_0)\begin{pmatrix}
v_1\\
v_2
\end{pmatrix},\begin{pmatrix}
v_1\\
v_2
\end{pmatrix}\rangle >\\
&\hspace{1.5cm} \delta_{x_0}(g(x_0)) \langle \begin{pmatrix}
\parcial{}{u_1}(g)^2&\parcial{}{u_1}(g)\parcial{}{u_2}(g)\\
\parcial{}{u_1}(g)\parcial{}{u_2}(g)&\parcial{}{u_2}(g)^2
\end{pmatrix}(x_0)\begin{pmatrix}
v_1\\
v_2
\end{pmatrix},\begin{pmatrix}
v_1\\
v_2
\end{pmatrix}\rangle \geq 0
.
\end{split}
\end{equation*}
This proves $\tilde{g}''_{x_0}(x)$ is strictly positive in an open
ball  $U_{x_0}$ around $x_0$ and multiplying $\tilde{g}''_{x_0}$
by a constant, if it is necessary, we have $\tilde{g}''_{x_0}(x)
\geq g''^{(i)}(x)$, for all $x \in U_{x_0}$. Since $\Gamma$ is
compact there exists a finite covering $\{U_{x_i}\}_{i=1}^N$ of
$\Gamma$, with associated functions $\{\tilde{g}_{x_i}\}_{i=1}^N$.
Let us define $\tilde{g}:= \sum_{i=1}^N \tilde{g}_{x_i}$. To see
that $\tilde{g}$ satisfies the lemma, it is enough to prove it for
$x$ in the set $A:=\{(y,u_1,u_2) \in Y \times
\R_+^2:k_1(3\epsilon,2\epsilon^2)u_2 \leq u_1 \leq
k_1(2\epsilon,3\epsilon^2)u_2\}\}$. Observe that for $x \in
\Gamma$, it follows by construction of $\tilde{g}$. Let
$(y,u_1,u_2) \in A$, then there exists $\lambda \in (0, \infty)$,
such that $(y,\lambda u_1,\lambda u_2) \in \Gamma$. Then, by
homogeneity, $(g^{(1)})''((y,u_1,u_2))=1/\lambda
(g^{(1)})''((y,\lambda u_1,\lambda u_2)) \leq 1/\lambda
\tilde{g}''(y,\lambda u_1,\lambda u_2)=\tilde{g}''(y,u_1,u_2)$.
$\square$
\subsection{Propagation observables}
Let $g$ be a Yafaev function. All our propagation observables are
derived from the following scaling of $g$, defined for $t>0$ and
$0<\delta<1$,
\begin{equation*}\label{eq:def gt}
g_t(x):=\begin{cases}
t^\delta g(y,t^{-\delta}u_1,t^{-\delta}u_2)&x=(y,u_1,u_2) \in  Y \times \R_+^2.\\
t^\delta g(z_i,t^{-\delta}u_i)&x=(z_i,u_i)\in Z_{i,0} \times
\R_+.\\
t^\delta \int \varepsilon_0
\varphi_0(\varepsilon_0)d\varepsilon_0& x\in X_0.
\end{cases}
\end{equation*}
We will be more precise about the value of $\delta$ later on. The
next results about  the derivatives of $g_t$ are the basis of
forthcoming estimates of propagation observables.
\begin{lem}{\upshape (cf. \cite[equation 7.18]{HS1})} \label{lem: propagation estimate g_t}
For each $((k_1,k_2),l) \in \N^2 \times \N$ and $t>0$  large
enough, there exist $C_1>0$ and $C_2>0$ such that:
\begin{equation*}\label{eq. estimate git}
\begin{split}
&\kappa(u_1)\kappa(u_2)\parcial{^{k_1}}{u_1^{k_1}}\parcial{^{k_2}}{u_2^{k_2}}
g^ {(j)}_t(x)\leq C_1 t^{\delta(1-\normv{k})}\text{ and }\\
& \hspace{2cm} \parcial{^l}{t^l} g^ {(j)}_t(x) \leq C_2
t^{\delta-l},
\end{split}
\end{equation*}
for $j=1,2,3$, for all $x \in X$ and $k_1 \geq 1$ or $k_2 \geq 1$.
\end{lem}
{\bf Proof:}
\\
We prove the lemma for $g^{(1)}$. The functions $g^{(2)}$ and
$g^{(3)}$ are treated in a similar form. Observe that the
integrand of (\ref{eq: def g1 int of epsilon}) has support in
$[2\epsilon,3\epsilon]$. Using Lebesgue dominated convergence
theorem, it is easy to see that there exists a $C>0$ depending
only on $k_1$ and $k_2$ such that:
\begin{equation}\label{eq:first ineq partial gt}
\begin{split}
&\normv{\parcial{^{k_2}}{u_2^{k_2}}\parcial{^{k_1}}{u_1^{k_1}}(g^{(1)})}\leq C \sum_{\normv{(j,s)}=k_1+k_2}\vert \int_{2\epsilon}^{3 \epsilon} \parcial{^{j_1}}{u_1^{j_1}}\left((1+\varepsilon_1)u_1\right) \varphi_1(\varepsilon_1)\\
&\hspace{0.5cm}\cdot  \parcial{^{j_0}}{u_1^{j_0}}\left( \Phi_0 \left((1+\varepsilon_1)t^{-\delta}u_1)\right) \right)\parcial{^{s_1}}{u_2^{s_1}}\parcial{^{j_2}}{u_1^{j_2}}\left(\Phi_2\left((1+\varepsilon_1)\frac{u_1}{u_2}-1)\right)\right)\\
&\hspace{0.5cm}\cdot
\parcial{^{s_2}}{u_2^{s_2}}\parcial{^{j_3}}{u_1^{j_3}}\left(
\Phi_3\left((1+\varepsilon_1)\frac{u_1}{\normv{u}}-1)\right)\right)d\varepsilon_1
\vert.
\end{split}
\end{equation}
We notice that the sum on the right-hand side of the above
inequality runs over the finite set of multi-indexes $(j,s)\in
\N^3 \times \N^2$ such that $\normv{(j,s)}= k_1+k_2$, where
$\normv{(j,s)}:=j_0+j_1+j_2+s_1+s_2$. We will denote by $B_{j,s}$
the terms of that sum and we will show that  they are uniformly
bounded by $t^ {\delta(1-k_1-k_2)}$. Since $g^{(1)}(z_1,u_1)=0$
for $ u_1 \leq \frac{l_0 t^{\delta}}{1+3\epsilon}$, the term
$B_{j,s}(y,u_1,u_2)=0$. Out of $k_1(3\epsilon,2\epsilon^2)u_2\leq
u_1 \leq k_1(2\epsilon,3\epsilon^2)u_2$ and $u_1 \geq
\frac{l_0t^\delta}{1+3\epsilon}$, the function $g^{(1)}$ is
constant or linear and the lemma follows easily. Hence we estimate
the terms $B_{j,s}$ only for $(y,u_1,u_2) \in Y \times \R_+\times
\R_+$ such that $k_1(3\epsilon,2\epsilon^2)u_2\leq u_1 \leq
k_1(2\epsilon,3\epsilon^2)u_2$ and $u_1 \geq
\frac{l_0t^\delta}{1+3\epsilon}$.
\\
\\
A direct computation shows that there exists a constant $C(j_0)$
such that:
\begin{equation} \label{eq:estimate gt Phi 0}
\begin{split}
\parcial{^{j_0}}{u_1^{j_0}}\left( \Phi_0 \left((1+\varepsilon_1)t^{-\delta}u_1)\right) \right)&\leq C(j_0)t^{-j_0 \delta }.
\end{split}
\end{equation}
We use above that $\varphi_0$ has compact support and hence all
its derivatives are bounded in $\R$. Observe that taking
$h(u_1,u_2):=(1+\varepsilon_1)\frac{u_1}{u_2}-1$ and
$f(v):=\frac{d^{j_2}}{dv^{j_2}}(\varphi_2)(v)$, one obtains:
\begin{equation*}\label{eq: estimating Phi2}
\frac{(1+\varepsilon_1)^{j_2}}{u_2^{j_2}}f \circ
h(u_1,u_2)=\parcial{^{j_2}}{u_1^{j_2}}\left(\Phi_2\left((1+\varepsilon_1)\frac{u_1}{u_2}-1)\right)\right).
\end{equation*}
Let $l\in \N$, let us define $M_{l}:=\{(k_1,\cdots,k_l)\in
\N^l:\sum_{i=1}^l ik_i=l\}$. We can conclude from  Fa\`a di
Bruno's formula  that for all $l \in \N$ and $\alpha \in M_l$
there exist constants $a_{l,k, \alpha}\in \R$ and  $C>0$ such that
\begin{equation} \label{eq: f circ g 1}
\begin{split}
&\normv{\parcial{^{s_1}}{u_2^{s_1}}\left(\frac{1}{u_2^{j_2}}f \circ h\right)(u_1,u_2)}  \leq C \sum_{l=0}^{s_1} \normv{\parcial{^l}{u_2^l} \left(f \circ h\right)(u_1,u_2) \frac{1}{u_2^{j_2+s_1-l}}}\\
&\leq C \sum_{l=0}^{s_1}\sum_{k=0}^{l}\sum_{\alpha \in M_l } \normv{a_{l,k,\alpha} (\partial^k(f)\circ h) (u_1,u_2) \Pi_{i=0}^{l}(\parcial{^i}{u_2^i}(h))^{\alpha_{i}}(u_1,u_2)\frac{1}{u_2^{j_2+s_1-l}}}\\
&\hspace{1.5cm}\leq  C \sum_{l=0}^{s_1}\sum_{k=0}^{l}\sum_{\alpha
\in M_l} \normv{a_{l,k,\alpha} (\partial^k(f)\circ h) (u_1,u_2)
\Pi_{i=0}^{l}
\frac{u_1^{\alpha_i}}{u_2^{(i+1)\alpha_i}}\frac{1}{u_2^{j_2+s_1-l}}}.
\end{split}
\end{equation}
For $u_1 \geq \frac{l_0 t^{\delta}}{1+3\epsilon}$ and  $u_1 \leq
k_1(2\epsilon,3\epsilon^2)u_2$,  there exists a constant
$C(s_1,j_2)>0$ such that the last term of (\ref{eq: f circ g 1})
is lower or equal than
\begin{equation}\label{eq:estimate gt Phi 2}
\begin{split}
 C \sum_{l=0}^{s_1}\sum_{k=0}^{l}\sum_{\alpha \in \N^n} (\partial^k(f)\circ h) (u_1,u_2) \Pi_{i=0}^{l} \frac{1}{u_2^{i\alpha_i+j_2+s_1-l}}
\leq C(s_1,j_2) t^{-\delta (j_2+s_1)},
\end{split}
\end{equation}
where we obtain the last inequality, since $j_2+s_1-l+\sum_{i=0}^l
i\alpha_i =j_2+s_1$ because the vectors $(\alpha_i) \in M_l$ and
the functions $\partial^k(f)$ have compact support.
\\
\\
Similar estimates can be done to obtain
\begin{equation} \label{eq:estimate gt Phi 3}
\begin{split}
\normv{\parcial{^{s_2}}{u_2^{s_2}}\parcial{^{j_3}}{u_1^{j_3}}\left(
\Phi_3\left((1+\varepsilon_1)\frac{u_1}{\normv{u}}-1)\right)\right)}
\leq Ct^{-\delta (j_3+s_2)}
\end{split}
\end{equation}
(\ref{eq:estimate gt Phi 0}), (\ref{eq: f circ g 1}),
(\ref{eq:estimate gt Phi 2}) and (\ref{eq:estimate gt Phi 3})
together with  (\ref{eq:first ineq partial gt}) imply the first
estimate of the lemma for the function $g^{(1)}$.
\\
\\
Next we will prove the second estimate of the lemma for the
function $g^{(1)}$. Let $\N \ni j \geq 1$, we proceed by induction
in $j$. The basis case, $j=1$, follows easily deriving with
respect to $t$ the scaling of expression (\ref{eq: def g1 int of
epsilon}). For $j \geq 1$, one uses Fa\`a di Bruno's formula for
$f=\varphi_0$ and $g(v)=(1+\varepsilon_1)t^{-\delta}v$, in a
similar way as it was used in (\ref{eq: f circ g 1}). One can
adapt the  proof  of the lemma for $g^{(1)}$ to the functions
$g^{(2)}$ and $g^{(3)}$. $\square$

We define the Heisenberg derivative of a function $h \in
C^\infty(\R_+\times X)$ by
\begin{equation}\label{eq: def Heissemberg derivative}
D_t h:=i[H,h]+\parcial{}{t} h.
\end{equation}
Now we estimate the first Heisenberg derivative $\gamma_t$ of
$g_t$ i.e. $$ \gamma_t:=D_t g_t=i[H,g_t]+\parcial{}{t} g_t.$$ We
will denote $\mathscr{W}_1(X,E)$ the domain of the self--adjoint
operator $\normv{H}^{1/2}$. Using an interpolation argument one
can see that $\mathscr{W}_1(X,E)$ coincides with the first Sobolev
space (see~\cite[Chapter 2]{Taylor}). Let us define the first
order differential operator $p:=i\begin{pmatrix}
\parcial{}{u_1}\\
\parcial{}{u_2}
\end{pmatrix}$ acting on sections $f \in C^\infty(Y\times \R_+\times \R_+,S)$
by $pf:=i\begin{pmatrix}
\parcial{}{u_1}f\\
\parcial{}{u_2}f
\end{pmatrix}$. We will denote $p^T$ the operator $i(
\parcial{}{u_1},
\parcial{}{u_2})$. The next lemma
shows that the asymptotic behavior of $\gamma_t$ is described by
the matrix function $g''_t(x)$ defined in (\ref{eq:defin
Hessian}), it is a consequence of lemma~\ref{lem: propagation
estimate g_t}.
\begin{lem}{\upshape (cf. \cite[equation (7.22)]{HS1})}\label{lem: estimate gammat-2partial g_t}
For all $2>\delta>0$ and all $\psi \in \mathscr{W}_1(X,E)$
\begin{equation*}
\langle D_t(\gamma_t-2\parcial{}{t} g_t) \psi_t,
\psi_t\rangle_{L^2(X,E)} =\langle \left( -4p^Tg''_tp+O(t^{-3
\delta})+O(t^ {\delta-2}) \right)\psi_t, \psi_t\rangle_{L^2(X,E)}
\end{equation*}
\end{lem}
{\bf Proof:}
\\
Observe that $\parcial{}{t}[H,g_t]=[H,\parcial{}{t}g_t]$, hence
$D_t(\gamma_t-2\parcial{}{t}
g_t)=-[H,[H,g_t]]-\parcial{^2}{t^2}g_t$. Using Leibnitz rule for
Laplacians and straightforward computations $$
[H,[H,g_t]]=4p^Tg''_t(x)p+\sum_{i,j=1}^2\partial_{jjii}(g_t).$$
According to lemma~\ref{lem: propagation estimate g_t},
$\partial_{jjii}(g_t)=O(t^{-3\delta})$ and $\parcial{^2}{t^2}g_t
\leq t^{\delta-2}$, which implies the lemma. $\square$

The next lemma is consequence of lemma~\ref{lem: estimate
gammat-2partial g_t}.
\begin{lem}{ \upshape (cf. \cite[theorem 7.5]{HS1})}\label{thm: int pgt'' is integrable}
For  $1>\delta>1/3$ there exists $C>0$  such that
\begin{equation*}
\normv{\int_1^\infty \langle p^Tg_t''p
\psi_t,\psi_t\rangle_{L^2(X,E)}dt} \leq C \norm{\psi}^2_1,
\end{equation*}
for all $\psi \in \mathscr{W}_1(X,E)$.
\end{lem}
{\bf Proof:}
\\
Using lemma~\ref{lem: estimate gammat-2partial g_t} we show
\begin{equation*}\label{eq: est gammat with g''}
\begin{split}
&\normv{\int_1^\infty \langle p^Tg_t''p
\psi_t,\psi_t\rangle_{L^2(X,E)} dt}\\
&\hspace{2cm}\leq \normv{\int_1^\infty \langle
D_t(\gamma_t-2\parcial{}{t} g_t) \psi_t,\psi_t \rangle_{L^2(X,E)}
dt}+K\norm{\psi}^2_{L^2(X,E)}
\end{split}
\end{equation*}
where $K>0$ is a constant. Next we estimate the first term in the
right side of the above inequality,
\begin{equation*}
\begin{split}
&\normv{\int_{1}^{t_0} \langle D_t(\gamma_t-2\parcial{}{t} g_t) \psi_t,\psi_t \rangle_{L^2(X,E)} dt}=\normv{\langle (\gamma_{t}-2\partial_{t} g_{t}) \psi_t,\psi_t \rangle_{L^2(X,E)} \vert^{t_0}_{t=1}\text{  }}\\
&\hspace{2cm}\leq \norm{(\gamma_{t}-2\partial_{t} g_{t})
\psi_t}_{L^2(X,E)}\text{  }\vert^{t_0}_{t=1} \cdot
\norm{\psi}_{L^2(X,E)} \leq C \norm{\psi}^2_1,
\end{split}
\end{equation*}
where the last inequality is true because lemma~\ref{lem:
propagation estimate g_t} implies that the first order
differential operator $\gamma_{t}-2\parcial{}{t} g_{t}$ has
bounded coefficients for $t \in [1,\infty)$ and hence it is
continuous from $L^2(X,E)$ to $\mathscr{W}_1(X,E)$. Since the
above inequality is true for arbitrary $t_0$ we have proved the
lemma.$\square$

We introduce and recall some  notation
\begin{equation*}
\begin{split}
g_{i,t}(x):=t^\delta g^{(i)}(t^{-\delta}x), \hspace{0.7cm}
g_t:=\sum_{i=0}^3g_{i,t}, \hspace{0.7cm} \gamma_{i,t}=D_t
g^{(i),t}, \hspace{0.7cm}\gamma_t:=\sum_{i=0}^3\gamma_{i,t},
\end{split}
\end{equation*}
where $D_t$ denotes the Heisenberg derivative  defined
in~(\ref{eq: def Heissemberg derivative}).

From part 3) of lemma~\ref{thm: properties of g}, it is easy to
see that $g''_t(x)$ is a positive matrix for all $t \in
[1,\infty)$ and $x \in X$. Therefore the matrix
$B(x,t):=\sqrt{g''_t(x)}$ is well defined. It is straightforward
to prove:
\begin{prop} \label{prop: pg p is equal pBp}
For $\varphi, \psi \in \mathscr{W}_1(X,E)$, the following equality
holds:
\begin{equation*}
\int_X \langle p^T g''_t p \psi, \varphi \rangle (x)
dvol(x)=\int_X \langle B p \psi,B p\varphi \rangle (x) dvol(x).
\end{equation*}
\end{prop}
Let  $\mathrm{Dom}(r)$ be the maximal domain in $L^2(X,E)$ of the
operator defined by multiplication by the function $r$ defined at
the beginning of section~~\ref{Mourre estimate}.
\begin{prop}
The domain  $\mathscr{W}_1(X,E) \cap \mathrm{Dom}(r)$ is invariant
under the action of $e^{iHt}$.
\end{prop}
{\bf Proof:}
\\
Let $\varphi \in \mathscr{W}_1(X,E) \cap \mathrm{Dom}(r)$. Since
$e^{iHt}$ and $H^{1/2}$ commute,  $e^{iHt} \varphi \in
\mathscr{W}_1(X,E)$, for all $t \in \R$. We have to show $r
e^{iHt} \varphi \in L^2(X,E)$. Let $\chi_n \in C^\infty_c(X)$ be
such that $\chi_n(x)=1$ for $x \in X_n$, and such that its
gradient $\nabla(\chi_n)$ and Laplacian $\Delta(\chi_n)$  are
bounded uniformly. We have
\begin{equation*} \label{eq: Dom r2 H1 is invariant 1}
\begin{split}
&\int_{X}\langle e^{iHt}\chi_nr^2 e^{-iHt}\varphi, \varphi \rangle (x)dvol(x)=\\
&i\int_X \int_0^t \langle e^{iHs}[H,\chi_n r^2]e^{-iHs}
\varphi,\varphi \rangle (x)dsdvol(x)+\int_X \chi_n r^2 \langle
\varphi,\varphi \rangle (x)dvol(x).
\end{split}
\end{equation*}
Let us see that the last integral is finite. By hypothesis $r
\varphi \in L^2(X,E)$, hence we can apply Lebesgue convergence
theorem to obtain
\begin{equation*}
\begin{split}
i\lim_{n \to \infty}\int_X \chi_n r^2 \langle \varphi,\varphi
\rangle(x) dvol(x)=\int_X r^2 \langle \varphi,\varphi \rangle(x)
dvol(x)<\infty.
\end{split}
\end{equation*}
Using that $[H,\chi_n r^2]$ is a first order differential operator
with uniformly bounded coefficients and Fubini's theorem we can
prove
$$i\int_X \int_0^t \langle e^{iHs}[H,\chi_n r^2]e^{-iHs}
\varphi,\varphi \rangle (x) dsdvol(x) \leq Ct \norm{\varphi}_1.$$
Lebesgue convergence theorem implies
\begin{equation*}
\begin{split}
&i\int_X \int_0^t \langle e^{iHs}[H, r^2]e^{-iHs}
\varphi,\varphi \rangle(x) dsdvol(x)\\
 &\hspace{1.5cm}=\lim_{n \to \infty} i\int_X \int_0^t \langle e^{iHs}
 [H,\chi_n r^2]e^{-iHs}
\varphi,\varphi \rangle(x) dsdvol(x)<\infty.\square\end{split}
\end{equation*}
The above proposition shows that the Heisenberg observables
$e^{iHt}\gamma_te^{-iHt}$ and $e^{iHt}g_te^{-iHt}$ are defined in
the dense domain $\mathscr{W}_1(X,E) \cap \mathrm{Dom}(r)$.
\begin{thm}{\upshape (cf. \cite[theorem 7.6]{HS1})} \label{thm: existence propagation observables}
1) The strong limits
\begin{equation*}
\gamma^+:=s-\lim_{t \to \infty} e^{iHt}\gamma_te^{-iHt},
\hspace{0.7cm } \gamma_k^+:=s-\lim_{t \to \infty}
e^{iHt}\gamma_{k,t}e^{-iHt}
\end{equation*}
exist on $\mathscr{W}_1(X,E)$  with respect to $L^2$-norm.

2) $\gamma^+$ and $\gamma_k^+$ have the following properties
\begin{equation*} \label{eq:gamma in terms g(t)}
\begin{split}
&\gamma^+_0=[\gamma^+,H]=[\gamma^+_k,H]=0,\hspace{0.4cm}\gamma^+
=s-\lim_{t \to \infty}
\frac{e^{iHt}g_te^{-iHt}}{t} \geq 0, \\
&\hspace{2cm}\gamma^+_k =s-\lim_{t \to \infty}
\frac{e^{iHt}g_{k,t}e^{-iHt}}{t} \geq 0 \text{ and }
\gamma^+=\sum_k \gamma^+_k.
\end{split}
\end{equation*}
where the last strong limits are taken over
$\mathscr{W}_1(X,E)\cap \mathrm{Dom}(r)$ with respect to the norm
$\norm{ \cdot}_{L^2(X,E)}$.

3) $\gamma^+$ and $\gamma_k^+$ are independent of $\delta\in
(1/3,1)$. Moreover, we have:
\begin{equation*} \label{eq: gamma+ does not depend of delta}
\gamma^+=s-\lim_{t \to \infty}e^{iHt}\frac{g(x)}{t} e^{-iHt},
\end{equation*}
where the strong limit is taken over $\mathscr{W}_1(X,E)\cap
\mathrm{Dom}(r)$, and where $g(x)$ is the unscaled Graf-Yafaev
function (similar roles play the functions $g^{(k)}$ for the
operators $\gamma_k^+$).
\end{thm}
Theorem~\ref{thm: existence propagation observables} will be
proved later on. We observe for the moment that from property 2)
we can deduce $\gamma^+_0=0$. Intuitively the importance of the
operators $\gamma^+_1,\gamma^+_2$ and $\gamma^+_3$ is that they
allow us to localize the absolutely continuous states of $H$ into
the regions
 $Z_1 \times \R_+$, $Z_2 \times \R_+$ and $Y \times \R^2_+$
associated with the domains of the operators $H_1,H_2$ and $H_3$.

We will use the following proposition to prove the existence of
$\gamma^+$.
\begin{prop} \label{In the strong limit gamma_t commutes with H} If one of the following limits exists,
 then $s-\lim_{t \to
\infty}e^{iHt} \gamma_t e^{-iHt}(H-\lambda)^{-2}=s-\lim_{t \to
\infty}(H-\lambda)^{-1}e^{iHt} \gamma_t e^{-iHt}(H-\lambda)^{-1}$.
\end{prop}
{\bf Proof:}
\\
We have that
\begin{equation*} \label{eq: prove exitence gamma 2}
\begin{split}
&(H-\lambda)^{-1}e^{iHt} \gamma_t
e^{-iHt}(H-\lambda)^{-1}=e^{iHt}(H-\lambda)^{-1}
\gamma_t(H-\lambda)^{-1}e^{-iHt}\\
&\hspace{0.8cm}=e^{iHt} \gamma_t (H-\lambda)^{-2}
e^{-iHt}-e^{iHt}(H-\lambda)^{-1}
[\gamma_t,H]e^{-iHt}(H-\lambda)^{-2}.
\end{split}
\end{equation*}
Then to prove the proposition, it is enough to prove:
\begin{equation} \label{eq: gamma + commutes with H}
\begin{split}
s-\lim_{t \to \infty}e^{iHt}(H-\lambda)^{-1}
[\gamma_t,H]e^{-iHt}(H-\lambda)^{-2}=0.
\end{split}
\end{equation}
By  lemma~\ref{lem: propagation estimate g_t},
$\norm{\parcial{}{t} (g_t)}_{0,0}=O(t^{\delta-1})$, where
$\norm{\cdot}_{0,0}$ denotes the norm of the bounded linear
operators acting in $L^2(X,E)$. Then we have:
\begin{equation*}\label{eq: prove exitence gamma 3}
\begin{split}
&s-\lim_{t \to \infty}e^{iHt}(H-\lambda)^{-1} [\gamma_t,H]e^{-iHt}(H-\lambda)^{-2}\\
&\hspace{0.7cm}=s-\lim_{t \to \infty}e^{iHt}
[(H-\lambda)^{-1},\gamma_t-\parcial{}{t}(g_t)]e^{-iHt}(H-\lambda)^{-1}.
\end{split}
\end{equation*}
Let $\psi:=(H-\lambda)^{-1}\varphi$, for $\varphi \in L^2(X,E)$.
We have
\begin{equation*}\label{eq:estimating H4 enough}
\begin{split}
&\norm{[(H-\lambda)^{-1},\gamma_t-\parcial{}{t}(g_t)]\psi}=\norm{[(H-\lambda)^{-1},[H,g_t]]\psi}_{L^2(X,E)}\\
&\hspace{0.6cm}\leq \norm{(H-\lambda)^{-1}[H,\sum_{i=1}^2\{-\partial_{ii}(g_t)-2\partial_i(g_t)\partial_i\}](H-\lambda)^{-1}\psi}_{L^2(X,E)}\\
&\hspace{0.6cm}\leq \sum_{j,i=1}^2
\norm{(H-\lambda)^{-1}\{\partial_{jjii}(g_t)+2\partial_{ji}(g_t)\partial_{ij}\}(H-\lambda)^{-1}\psi}_{L^2(X,E)}.
\end{split}
\end{equation*}
Using lemma \ref{lem: propagation estimate g_t}, one can prove
$\norm{\partial_{iijj}(g_t)}_{0,0} \leq C t^{-3\delta}$, that
implies:
\begin{equation*}\label{eq: estimate partialjjii}
\begin{split}
\norm{(H-\lambda)^{-1}\partial_{jjii}(g_t)(H-\lambda)^{-1}\psi}_{L^2(X,E)}\leq
C t^{-3 \delta}.
\end{split}
\end{equation*}
Now we analyze the term
$\norm{(H-\lambda)^{-1}\partial_{ji}(g_t)\partial_{ij}(H-\lambda)^{-1}\psi}_{L^2(X,E)}$.
Since $\partial_{ji}(g_t)\partial_{ij}$ is a second order
differential operator with coefficients bounded uniformly in $x
\in X$ and $t \in \R$,  it defines a continuous operator from
$\mathscr{W}_2(X,E)$ to $L^2(X,E)$, hence:
\begin{equation*}
\begin{split}
&\norm{(H-\lambda)^{-1}\partial_{ji}(g_t)\partial_{ij}(H-\lambda)^{-1}\psi}_{L^2(X,E)} \\
&\leq  \norm{(H-\lambda)^{-1}}_{0,2}\cdot
\norm{\partial_{ij}(g_t)\partial_{ij}}_{2,0}\cdot
\norm{(H-\lambda)^{-1}}_{0,2}\cdot \norm{\psi}_{L^2(X,E)},
\end{split}
\end{equation*}
where $\norm{\cdot}_{k,l}$ denotes the operator norm from
$\mathscr{W}_{k}(X,E)$ to $\mathscr{W}_{l}(X,E)$. We observe that,
by lemma \ref{lem: propagation estimate g_t}, we have $
\norm{\partial_{ij}(g_t)\partial_{ij}}_{2,0} \leq C t^{\delta-1}$,
this finishes the proof of the proposition. $\square$

{\bf Proof of theorem~\ref{thm: existence propagation
observables}:}
\\
\\
{\bf 1. Existence of $\gamma^+$ and  $\gamma_k^+$:}
Lemma~\ref{lem: propagation estimate g_t} implies that
$([H,g_t])_{t \in \R_+}$ and $(\parcial{}{t}(g_t))_{t \in \R_+}$
have coefficients bounded uniformly in $t \in \R_+$ and $x \in X$
 and then we can
deduce the inequalities
\begin{equation*}
\begin{split}
\norm{e^{iHt} \parcial{}{t}(g_t) e^{-iHt}\varphi}_{L^2(X,E)} \leq
C
t^{\delta -1} \norm{\varphi}_{L^2(X,E)}\\
\norm{e^{iHt}[H,g_t] e^{-iHt}\varphi}_{L^2(X,E)} \leq C\norm
{\varphi}_1,
\end{split}
\end{equation*}
for all $\varphi \in \mathscr{W}_4(X,E)\subset\mathscr{W}_1(X,E)$.
The previous estimates show that, assuming the existence of the
limit, the function $\varphi \mapsto \lim_{t \to \infty} e^{iHt}
\gamma_t e^{-iHt}\varphi$ would be a continuous linear map (as a
function from $\mathscr{W}_1(X,E)$ to $L^2(X,E)$). Since
$\mathscr{W}_2(X,E)\subset\mathscr{W}_1(X,E)$ is dense with
respect to the first Sobolev norm $\norm{\cdot}_1$, it is enough
to prove that the limit $\lim_{t\to \infty}e^{iHt} \gamma_t
e^{-iHt}(H-\lambda)^{-2}$ exists and hence, by proposition~\ref{In
the strong limit gamma_t commutes with H}, it is enough to prove
the existence of the limit $s-\lim_{t \to
\infty}(H-\lambda)^{-1}e^{iHt}\gamma_t e^{-iHt}(H-\lambda)^{-1}$
with respect to the norm $\norm{\cdot }_{L^2(X,E)}$. Since
$\norm{\parcial{}{t} g_t}_{0,0}=O(t^{\delta-1})$,  we have
\begin{equation*}\label{eq: gammat equiv gammat-partialt}
\begin{split}
s-\lim_{t \to \infty}&(H-\lambda)^{-1}e^{iHt}\gamma_t e^{-iHt}(H-\lambda)^{-1}=\\
&s-\lim_{t \to
\infty}(H-\lambda)^{-1}e^{iHt}(\gamma_t-2\parcial{}{t} (g_t))
e^{-iHt}(H-\lambda)^{-1}.
\end{split}
\end{equation*}
We will show the existence of the last limit with respect to the
$L^2$-norm. We  denote $\tilde{\gamma}_t:=\gamma_t-2\parcial{}{t}
(g_t)$.
\\
\\
Define $\varphi_t:=(H-\lambda)^{-1}e^{iHt}\tilde{\gamma}_t
e^{-iHt}(H-\lambda)^{-1}\psi$ for $\psi \in L^2(X,E)$. We will
prove  that $\int_1^\infty
\norm{\parcial{}{t}\varphi_t}_{L^2(X,E)} dt$  is finite. Observe
that
\begin{equation*}
\begin{split}
\parcial{}{t}\varphi_t:=(H-\lambda)^{-1}e^{iHt}D_t\tilde{\gamma}_t e^{-iHt}(H-\lambda)^{-1}\psi.
\end{split}
\end{equation*}
From lemma~\ref{lem: estimate gammat-2partial g_t}, for
$\delta>1/3$, we can deduce:
\begin{equation*}
\begin{split}
D_t\tilde{\gamma}_t=p^Tg''_tp+L^2\text{-norm integrable in $t$
terms}.
\end{split}
\end{equation*}
Therefore it remains to prove that $ u_t:=(H-\lambda)^{-1}e^{iHt}
p^Tg''_tp e^{-iHt}(H-\lambda)^{-1}\psi$ is $L^2$-norm integrable
in $[1,\infty)$. We use  Cauchy--Schwarz inequality and
proposition~\ref{prop: pg p is equal pBp} to prove
\begin{equation*} \label{eq: existence ut proof gamma+ existe}
\begin{split}
&\int_{1}^{s}\norm{u_t}_{L^2(X,E)}^2dt=\int_{1}^{s}\sup_{\norm{v}_{L^2(X,E)}=1} \normv{\langle v,u_t \rangle_{L^2(X,E)}}^2dt\\
&\hspace{1cm}\leq
\sup_{\norm{v}_{L^2(X,E)}=1}\int_{1}^{s}\norm{B_t p
e^{-iHt}(H-\overline{\lambda})^{-1}v}_{L^2(X,E)}^2dt\\
&\hspace{3cm}\cdot \int_{1}^{s}
\norm{B_tpe^{-iHt}(H-\lambda)^{-1}\psi}_{L^2(X,E)}^2 dt.
\end{split}
\end{equation*}
By lemma~\ref{thm: int pgt'' is integrable} the last two integrals
are  bounded; hence $u_t$ is $L^2$--norm integrable in $t$.  We
have proved the existence of $\gamma^+$, the existence of
$\gamma_k^+$ is proved following a very similar reasoning.
\\
\\
{\bf 2.  Proof of parts 2) and 3) of theorem~\ref{thm: existence
propagation observables}:} Since $\gamma^+$ exists on
$\mathscr{W}_1(X,E)$, it follows from (\ref{eq: gamma + commutes
with H}) that $\gamma^+
(H-\lambda)^{-1}=(H-\lambda)^{-1}\gamma^+$. Hence
$[\gamma^+,H]=(H+\lambda)\{(H+\lambda)^{-1}\gamma^+-\gamma^+(H+\lambda)^{-1}\}(H+\lambda)=0$.
A similar proof applies for $\gamma_k^+$.
\\
\\
Now we prove that  $\lim_{t \to
\infty}\frac{e^{iHt}g_te^{-iHt}}{t}\varphi=\gamma^+\varphi $ for
$\varphi \in \mathrm{Dom}(r)\cap \mathscr{W}_1(X,E)$ and  where
the limit is considered in the $L^2$-norm. Using that
$e^{iHt}\gamma_te^{-iHt}=\parcial{}{t}e^{iHt}g_te^{-iHt}$, we have
$$\gamma^+=s-\lim_{t \to \infty} \frac{1}{t} \int_1^t
\parcial{}{s}e^{iHs}g_se^{-iHs}ds=s-\lim_{t \to \infty} \frac{e^{iHt}g_te^{-iHt}}{t} \geq 0.$$
Finally we prove part 3) of theorem~\ref{thm: existence
propagation observables}. Observe that $g_t=g$ for $x \in X-X_R$
and for $R>\frac{l_0+\epsilon}{1+2\epsilon^2}$; hence
$t^{-1}\norm{g_t-g}_{L^2(X,E)} \leq C t^{\delta-1}$. Part 3)
follows from part 2) of the theorem and this fact. $\square$
\subsection{Propagation observables and Mourre's inequality}
Next we discuss the connection between the operator $\gamma^+ $
and  Mourre's inequality enunciated in theorem~\ref{thm:Mourre
estimate}.
\begin{defin}{\upshape (cf. \cite[(6.17)]{HS1})}\label{defin:Mourre interval}
A finite, open interval $I \subset \R$ will be called {\bf a
Mourre interval} if for all $\psi \in E_I(H) \cap \mathrm{Dom}(r)$
\begin{equation*}\label{eq:def: weak Mourre}
 \langle E_I(H)i[H,i[H,r^2]]E_I(H)\psi,\psi\rangle_{L^2(X,E)} \geq C \langle \psi,\psi\rangle_{L^2(X,E)} \text{ for some } C>0.
\end{equation*}
\end{defin}
\begin{lem}{\upshape (cf. \cite[lemma 7.7]{HS1})} \label{gamma+ is invertible in Mourre intervals}
Let $\mathscr{H}_I:=E_I(H)$ be the spectral subspace of $H$
associated to a  Mourre interval $I$. Then ${\gamma^+}^2$ reduces
to a strictly positive operator $\mathscr{H}_I \to \mathscr{H}_I$.
In particular $\mathscr{H}_I \subset \mathrm{Im}(\gamma^+)$.
\end{lem}
{\bf Proof:}
\\
According to theorem~\ref{thm: existence propagation observables},
$\gamma^+$ is $H$-bounded and commutes with $H$, then it reduces
to $\mathscr{H}_I \to \mathscr{H}_I$. Let $\psi \in
\mathscr{H}_I$, by theorem~\ref{thm: existence propagation
observables}
\begin{equation*}
\begin{split}
\langle \psi,&{\gamma^+}^2 \psi \rangle_{L^2(X,E)}=\lim_{t \to \infty}\frac{1}{t^2 } \langle e^{2iHt}\psi,g_t^{2}e^{2iHt}\psi \rangle_{L^2(X,E)}\\
&\hspace{2cm} \geq \lim_{t \to \infty}\frac{1}{t^2 }\langle
e^{2iHt}\psi,r^{2}e^{2iHt}\psi \rangle_{L^2(X,E)}.
\end{split}
\end{equation*}
Define the function $h(t):=\langle
e^{2iHt}\psi,r(x)^{2}e^{2iHt}\psi \rangle_{L^2(X,E)}$. Since $I$
is a Mourre interval, there exists $c>0$ such that $h''(t) \geq
c>0$. Then, there exist $c_1 \in \R$ and $c_2 \in \R$ such that
$h(t)\geq ct^2+c_1t +c_2$ and
$$\lim_{t \to \infty}\frac{1}{t^2
}\langle e^{2iHt}\psi,r(x)^{2}e^{2iHt}\psi\rangle_{L^2(X,E)} \geq
c>0. \square
$$

As a consequence of theorem~\ref{thm:Mourre estimate} we have that
if $\lambda \in \R$ is not an $L^2$--eigenvalue nor a threshold of
$H$, then $\lambda$ belongs to some Mourre interval $I$. In
\cite{CanoMourre} and \cite{CANOTHESIS} is proved by different
methods that the set of  $L^2$--eigenvalues of $H$ is countable
and it accumulates only in the set of thresholds
 $\sigma_{pp}(H^{(1)})\cup\sigma_{pp}(H^{(2)})\cup
\sigma_{pp}(H^{(3)})$. The next corollary follows from these
facts.
\begin{cor}{\upshape (cf. \cite[page 3480]{HS1})} \label{cor: Mourre intervals are dense L2ac}
The sum of eigenspaces $E_I(H)$, associated to  Mourre intervals
$I$, is a $L^2$-dense set on the absolutely continuous part of $H$
\end{cor}
\subsection{Deift-Simon wave operators}
The proof of the following theorem follows the same lines  of  the
proof of the existence of $\gamma^+$ and $\gamma_k^+$ in
theorem~\ref{thm: existence propagation observables} and the proof
of similar facts given in~\cite[page 3492]{HS1}, because of this
we omit the proof here.
\begin{thm}{\upshape (cf. \cite[page 3492]{HS1})}\label{thm: existence Deift}
For $k=1,2,3$, the Deift-Simon wave operators,
\begin{equation*}
\omega_k:=s-\lim_{t \to \infty} e^{iH_k t} \gamma_{k,t}e^{-iHt},
\end{equation*}
exist, with respect to the $L^2$-norm, on $\mathscr{W}_1(X,E)$ for
$\delta$ satisfying $\min(3 \delta, 2- \delta)<1$.
\end{thm}
As we explained below theorem~\ref{thm: existence propagation
observables}, intuitively the importance of the operators
$\gamma_{k,t}$ is that they allow us to localize in the domains of
the operators $H_k$ the absolutely continuous states of $H$. In
theorem~\ref{thm: existence Deift} we  find states whose dynamics
under $H_k$ behave asimptotically as the dynamic of these
localizations under $H$. We will formalize these intuitions in the
next section.
\subsection{Proof of asymptotic clustering}
\label{section: asymp clustering}In this section we prove
asymptotic clustering wich finishes the proof of
theorem~\ref{thm:asympt completeness}. We say that {\bf $\psi \in
L^2_{ac}(X,E)$ clusters asymptotically, } if there exist
$\varphi_k \in L^2_{pp}(Z_k,E_k) \otimes L^2(\R_+)$ for $k=1,2$
and $\varphi_3 \in L^2(Y \times \R_+^2,E)$ such that
equation~(\ref{eq: asym completeness}) holds.

Let $\psi \in E_I(H) \cap \mathscr{W}_2(X,E)$ for $I$ a Mourre
interval as defined in definition~\ref{defin:Mourre interval}. By
lemma~\ref{gamma+ is invertible in Mourre intervals} and
theorem~\ref{thm: existence propagation observables}, we have
\begin{equation*}
\begin{split}
\psi=\sum_{k=1}^3 \gamma_k^+ \varphi \approx \sum_{k=1}^3 e^{iHt}
\gamma_{k,t} e^{-iHt}\varphi,
\end{split}
\end{equation*}
where $\approx$ means that the difference of the two related
expressions vanishes in $L^2$-norm as $t \to \infty$.
Theorem~\ref{thm: existence Deift} implies
\begin{equation} \label{eq: decomposition using Deift1}
\begin{split}
\psi_t \approx \sum_{i=1}^3 e^{-iH_k t} \varphi_k, \text{ for  }
\varphi_k:= \omega_k^+ \varphi,
\end{split}
\end{equation}
that with corollary~\ref{cor: Mourre intervals are dense L2ac}
imply that the wave operators $W_\pm(H_k,H)$ exist, for $k=1,2,3$.
\begin{prop}\label{prop: asympt completeness eHkt version}For all $\psi \in L^2_{ac}(X,E)$  there exist  $\varphi_k \in L^2(Z_k \times \R_+,E)$,  for $k=1,2,3$,  such that
\begin{equation*}
\begin{split}
\lim_{t \to \infty} \norm{e^{\pm iHt}\psi-\sum_{k=1}^3e^{\pm
iH_kt}\varphi_k}_{L^2(X,E)}=0.
\end{split}
\end{equation*}
\end{prop}
Proposition~\ref{prop: asympt completeness eHkt version} is a kind
of asymptotic completeness, however the sum of the wave operators
$W_{\pm}(H,H_k)$ ($k=1,2,3$) is not a direct sum, since their
images are not necessarily orthogonal.

For $k=1,2$ and the $\varphi_k$'s of (\ref{eq: decomposition using
Deift1}), we have
$\varphi_k=\Pi_{k,pp}\varphi_k+\Pi_{k,ac}\varphi_k$, where
$\Pi_{k,pp}$ and $\Pi_{ac,d}$ denote the orthogonal projection
over the closed subspaces of $L^2(X,E)$, $L^2_{pp}(Z_k,E_k)
\otimes L^2(\R_+)$ and $L^2_{ac}(Z_k,E_k) \otimes L^2(\R_+)$. It
is easy to see that $e^ {\pm it H_k}\varphi_k=e^ {\pm it
H_{k,pp}}\Pi_{k,pp}\varphi_k+e^ {\pm it
H_{k,ac}}\Pi_{k,ac}\varphi_k$.

Since $\Pi_{k,ac}\varphi_k \in L^2_{ac}(Z_k,E_k) \otimes
L^2(\R_+)$ and $W_\pm(H_k,H_3)$ is an isometry, there exists
$\tilde{\varphi}_k \in L^2(Y \times \R_+^2,E)$ such that
$\Pi_{k,ac}\varphi_k=W_\pm(H_k,H_3) \tilde{\varphi}_k$. We
conclude
\begin{equation*} \label{eq: using WHk H3}
\begin{split}
e^{\pm iHt}\psi-\sum_{k=1}^3e^{\pm iH_kt}\varphi_k&=e^{\pm iHt}\psi-\sum_{k=1}^2\{e^{\pm iH_{k,ac} t}W_\pm(H_k,H_3) \tilde{\varphi}_k\\
&-e^{\pm iH_{k,pp} t}\Pi_{k,pp}\varphi_k\}-e^{i\pm H_3
t}\varphi_3.
\end{split}
\end{equation*}
Observe that
\begin{equation*}
\begin{split}
\lim_{t \to \infty}\norm{e^{\pm iH_{k,ac} t}W_\pm(H_k,H_3)
\tilde{\varphi}_k-e^{\pm iH_3t}\tilde{\varphi}_k}_{L^2(X,E)}=0,
\end{split}
\end{equation*}
for $k=1,2$. The above computations imply
\begin{prop}\label{prop: eiHt in term eiHkd and eiH3}
For all $\psi \in L^2_{ac}(X,E)$  there exist $\phi_k \in
L^2_{pp}(Z_k,E_k) \otimes L^2(\R_+)$, for $k=1,2$, and $\varphi
\in L^2(Y \times \R_+^2,E)$,  such that
\begin{equation*}
\begin{split}
\lim_{t \to \infty} \norm{e^{\pm iHt}\psi-e^{\pm
iH_3t}\varphi-\sum_{k=1}^2e^{\pm iH_{k,pp}}\phi_k}_{L^2(X,E)}=0.
\end{split}
\end{equation*}
\end{prop}
Let $\varphi \in L^2(Y \times \R_+^2,E)$, we can calculate
\begin{equation*}\label{eq:toward asym of eIH3 in terms WHkH3}
\begin{split}
e^{\pm iH_3t}\varphi=&e^{\pm iH_3t}\varphi- e^{\pm iH_{1,ac}}W_\pm(H_1,H_3)\varphi\\
&+e^{\pm iH_3t}\varphi- e^{\pm iH_{2,ac}}W_\pm(H_2,H_3)\varphi\\
&-e^{\pm iH_3t}\varphi+ \sum_{k=1}^2 e^{\pm
iH_{k,ac}}W_\pm(H_k,H_3)\varphi.
\end{split}
\end{equation*}
Observe that for all $\varphi \in L^2(Y \times \R_+^2,E)$, we have
\begin{equation*}\label{eq: eiH3 asym WHkH3}
\begin{split}
\lim_{t \to \infty} \norm{e^{\pm iH_3t}\varphi-e^{\pm
iH_{k,ac}}W_\pm(H_k,H_3)\varphi}_{L^2(X,E)}=0,
\end{split}
\end{equation*}
for $k=1,2$, hence,
\begin{equation*}\label{eq:eiH3 in terms eiHkac}
\begin{split}
e^{\pm iH_3t}\varphi \approx &-e^{\pm iH_3t}\varphi+ \sum_{k=1}^2
e^{\pm iH_{k,ac}}W_\pm(H_k,H_3)\varphi.
\end{split}
\end{equation*}
Proposition \ref{prop: eiHt in term eiHkd and eiH3} and the
previous computation imply asymptotic clustering and hence
theorem~\ref{thm:asympt completeness}.

Let us denote
$$\mathscr{W}_\pm:=W_{\pm}(H,H_3) \oplus \bigoplus_{k=1}^2
W_\pm(H,H_{k,pp}),$$ acting from $L^2(Y \times \R_+^2)\oplus
\bigoplus \left(L^2_{pp}(Z_k,E_k)\otimes L^2(\R_+) \right)$ to
$L^2_{ac}(X,E)$. We define the scattering operator
\begin{equation}\label{eq:defin scattering}
\mathscr{S}:=\left(\mathscr{W}_- \right)^{-1}\mathscr{W}_+
\end{equation}
In a forthcoming article, we plan to study  how the scattering
operator $\mathscr{S}$ encodes geometric information, particularly
we would like to generalize the approach of~\cite{MuellerCorner}
to prove a signature formula that would be closely related with
the formulas of~\cite{L2signMazzeo}.

\appendix 
\section{Stationary phase methods}\label{app: stationary phase}
Let $V \in C^\infty_c(\R)$. In this appendix we prove
$\int_{-\infty}^\infty \norm{Ve^{it\tilde{b}}u}dt < \infty$ where
$\tilde{b}$ is the self--adjoint operator associated to
$-\frac{d^2}{dx^2}:C^\infty_c(\R) \to L^2(\R)$. We use stationary
phase methods as explained in~\cite{RS3}. Let $u \in
\mathscr{S}(\R)$ be such that $\hat{u}$ has compact support
contained in an interval $[a,d]$. Here $\hat{u}$ denotes the
Fourier transform of $u$. Let
$$u_t(x):=e^{it \tilde{b}}u=(\frac{1}{2\pi})^{1/2} \int
\exp[it(xk-tk^2)]\hat{u}(k)dk.$$ From \cite[Corollary, page
38]{RS3} we see that for all $m$ there exists a $c$ depending on
$m$, $u$ and the interval $[a,b]$ such that
$$\normv{u_t(x)}\leq c(1+\normv{x}+t)^{-m}$$ for
all $x, t$ such that $x/t \notin [a,d]$. From this we deduce that
\begin{equation}\label{eq: est Vut out of at bt}
\begin{split}
&(\int_{-\infty}^{at} +\int_{dt}^{\infty})
\normv{V(x)}^2\normv{u_t(x)}^2dx\leq c  (1+t)^{-2}.\end{split}
\end{equation}
\cite[Corollary, page 41]{RS3} proves $\normv{u_{\pm t}(x)}^2 \leq
C t^{-1}$ for $t>1$, then $$\int_1^\infty \left( \int_{at}^{dt}
\normv{V(x)}^2 \normv{u_{t}(x)}^2 dx\right)dt \leq c \int_1^\infty
t^{-1/2}\left( \int_{at}^{dt} \normv{V(x)}^2 dx\right)dt.
$$
Making the change of variables $x=xt$ we obtain that for all $m
\in \N $  there exists a $C$ such that $\normv{\int_{at}^{dt}
\normv{V(x)}^2 dx} \leq \frac{Ct}{1+t^m}$, then
\begin{equation}\label{eq: est Vut  in at bt}
\begin{split}\int_{at}^{dt} \normv{V(x)}^2\normv{u_t(x)}^2dx\leq
  C \frac{t^{1/2}}{1+t^5}.
\end{split}
\end{equation}
(\ref{eq: est Vut out of at bt}) and (\ref{eq: est Vut  in at bt})
show that $\int_{-\infty}^\infty \norm{Ve^{it\tilde{b}}u}dt <
\infty$.
\\
\\
Next we make some classical comments in order to  extend the
previous estimates to $\norm{Ve^{itb}\varphi}$ where $b$ is the
self--adjoint operator associated to
$-\frac{d^2}{dx^2}:C^\infty_c(\R_+) \to L^2(\R_+)$ with Dirichlet
boundary conditions at $0$. We observe that for all $u \in
\mathscr{S}((0,\infty))$ such that $\hat{u} \in C^\infty((0,
\infty))$, the function
$$\tilde{u}(x):=\begin{cases}u(x)& x\in (0,
\infty),\\
0& x=0,\\
-u(-x)& \text{otherwise}
\end{cases}$$
is an odd function in $\mathscr{S}(\R)$ such that
$\hat{\tilde{u}}$ has compact support. Since $\tilde{u}$ is odd,
$\hat{\tilde{u}}=2i \int_0^\infty \sin(xy)u(y)dy$. From these
observations, (\ref{eq: est Vut out of at bt}) and (\ref{eq: est
Vut in at bt}), we deduce $\int_{-\infty}^\infty
\norm{Ve^{itb}u}dt < \infty$.

\bibliographystyle{alpha}
\bibliography{literature}
\end{document}